\documentclass{article}

\usepackage{latexsym, amsfonts, amsmath, amssymb}
\newcommand{\be}{\begin{equation}}
\newcommand{\ee}{\end{equation}}
\newcommand{\bea}{\begin{eqnarray}}
\newcommand{\eea}{\end{eqnarray}}
\newcommand{\bean}{\begin{eqnarray*}}
\newcommand{\eean}{\end{eqnarray*}}

\newcommand{\brray}{\begin{array}}
\newcommand{\erray}{\end{array}}
\newcommand{\ben}{\begin{equation}{nonumber}}
\newcommand{\een}{\end{equation}{nonumber}}

\newtheorem{dfn}{Definition}[section]
\newtheorem{thm}[dfn]{Theorem}
\newtheorem{lmma}[dfn]{Lemma}
\newtheorem{ppsn}[dfn]{Proposition}
\newtheorem{crlre}[dfn]{Corollary}
\newtheorem{xmpl}[dfn]{Example}
\newtheorem{rmrk}[dfn]{Remark}
\newcommand{\bdfn}{\begin{dfn}}
\newcommand{\bthm}{\begin{thm}}
\newcommand{\blmma}{\begin{lmma}}
\newcommand{\bppsn}{\begin{ppsn}}
\newcommand{\bcrlre}{\begin{crlre}}
\newcommand{\bxmpl}{\begin{xmpl}}
\newcommand{\brmrk}{\begin{rmrk}}
\newcommand{\edfn}{\end{dfn}}
\newcommand{\ethm}{\end{thm}}
\newcommand{\elmma}{\end{lmma}}
\newcommand{\eppsn}{\end{ppsn}}
\newcommand{\ecrlre}{\end{crlre}}
\newcommand{\exmpl}{\end{xmpl}}
\newcommand{\ermrk}{\end{rmrk}}

\newcommand{\IR}{\mathbb{R}}

\newcommand{\IZ}{\mathbb{Z}}

\newcommand{\cla}{{\cal A}}

\newcommand{\clf}{{\cal F}}

\newcommand{\clh}{{\cal H}}

\newcommand{\clq}{{\cal Q}}

\newcommand{\cls}{{\cal S}}

\newcommand{\clu}{{\cal U}}
\newcommand{\clv}{{\cal V}}

\def\a*{{\cal A}_{h,*}}
\def\B{{\cal B}(h)}
\def\B1{{\cal B}_1(h)}
\def\b{{\cal B}^{\rm s.a.}(h)}
\def\b1{{\cal B}^{\rm s.a.}_1(h)}

\newcommand{\ot}{\otimes}

\def \qed {$\Box$}

\def\a*{{\cal A}_{h,*}}
\def\B{{\cal B}(h)}
\def\B1{{\cal B}_1(h)}
\def\b{{\cal B}^{\rm s.a.}(h)}
\def\b1{{\cal B}^{\rm s.a.}_1(h)}

\begin{document}
\begin{center}
{\Large{\bf Quantum Isometry groups of dual of finitely generated discrete groups and quantum groups }}\\~\\ 
{\large {Debashish Goswami\footnote{Partially supported by Swarnajayanti fellowship from D.S.T (Govt of India)} and Arnab Mandal}}\\
Indian Statistical Institute\\
203, B. T. Road, Kolkata 700108\\
Email: goswamid@isical.ac.in\\
\end{center} 
\begin{abstract}
 We study quantum isometry groups, denoted by  $\mathbb{Q}(\Gamma, S)$,  of spectral triples on $C^*_r(\Gamma)$ for a finitely generated discrete group $\Gamma$ 
 coming from the word-length metric with respect to a symmetric generating set $S$.  We first prove a few general results about $\mathbb{Q}(\Gamma, S)$ 
  including:
  \begin{itemize}
 \item For a group $\Gamma$ with polynomial growth property, the dual of $\mathbb{Q}(\Gamma, S)$ has polynomial growth property provided the action of $\mathbb{Q}(\Gamma,S)$ on $C^*_r(\Gamma)$ has full spectrum.  
   \item $\mathbb{Q}(\Gamma, S) \cong QISO(\hat{\Gamma}, d)$ for any discrete abelian group $\Gamma$, where $d$ is a  suitable metric on the dual compact abelian group $\hat{\Gamma}$.
   \end{itemize}  
  We then carry out explicit computations  of $\mathbb{Q}(\Gamma,S)$ for several classes of 
   examples including free and direct product of cyclic groups, Baumslag-Solitar group, Coxeter groups etc. In particular, we have computed quantum isometry groups of  all finitely 
    generated abelian groups which do not have factors of the form $\IZ_2^k$ or $\IZ_4^l$ for some $k, l$ in the direct product decomposition into cyclic subgroups.

\end{abstract}

 \section{Introduction}
 It is a very interesting problem, both from the physical and mathematical viewpoint, to understand and classify quantum symmetries of possibly non-commutative $C^*$-algebras (usually with further structures), i.e. possible actions of quantum groups on them. In \cite{Manin}, this problem was considered in an algebraic and categorical setting, leading to the realization of some of the well known (algebraic) quantum groups such as $SL(2,q)$ as the universal object in some category of quantum groups acting on the quantum 2-plane. S.Wang \cite{wang} took up a similar problem in the analytical framework of compact quantum groups acting on $C^*$-algebras. Later on, a number of mathematicians including Wang, Banica, Bichon and others (\cite{wang}, \cite{ban_1}, \cite{finite graph}) developed a  theory of quantum automorphism groups of finite dimensional $C^*$-algebras as well as quantum isometry groups of finite metric spaces and finite graphs. In \cite{Gos} the first named author of the present article extended 
such 
constructions to the set up of possibly infinite dimensional $C^*$-algebras, and more interestingly, that of spectral triples a la Connes \cite{con}, by defining and studying quantum isometry groups of spectral triples. This led to the study of such quantum isometry groups by many authors including Goswami, Bhowmick, Skalski, Banica, Bichon, Soltan, Das, Joardar and others.
 In the present article, our focus is on a rather special yet interesting and important class of spectral triples, namely 
 those coming from the word-length metric of finitely generated discrete groups with respect to some given symmetric generating set. There have been several articles already on computations and study
 of the quantum isometry groups of such spectral triples, e.g \cite{grp algebra}, \cite{dihedral}, \cite{S_n}, \cite{free cyclic}, \cite{org filt} and references therein. 
 However, for a systematic and unified study of quantum isometry groups of such spectral triples, one needs to look at many more examples and then try to identify some general pattern. 
  This is the main objective of this paper. We have not yet been able to propose a general theory, but almost complete the understanding of quantum isometry groups of direct and free 
   product of cyclic groups, except a few cases only. Besides, we treat several other important classes of groups. 
 
 We begin by proving some general facts about the quantum isometry group $\mathbb{Q}(\Gamma,S)$ of a discrete group $\Gamma$ with a finite symmetric generating set $S$. We prove, among other things,
  the following two interesting results:
  \begin{enumerate}
 \item If $\Gamma$ has polynomial growth and the action of $\mathbb{Q}(\Gamma, S)$ on $C^*_r(\Gamma)$ has full spectrum, then the dual (discrete quantum group)  of 
   $\mathbb{Q}(\Gamma)$ has polynomial growth.
   \item In case $\Gamma$ is abelian, there is a metric on the dual compact abelian group $\hat{\Gamma}$ such that the corresponding quantum isometry group in the metric space sense 
    (as in \cite{metric iso}) exists and coincides with $\mathbb{Q}(\Gamma,S)$.
    \end{enumerate} 
 
 Next we carry out several explicit computations. We have given special emphasis on groups of the 
 form $\Gamma_1 \ast \Gamma_2 \ast \cdot\cdot \ast \ \Gamma_k$ or $\Gamma_1 \times \Gamma_2 \times \cdot\cdot \times \ \Gamma_k$, 
 where $\Gamma_i= \mathbb{Z}_{n_i}$ for some $n_i$. We have proved that in many cases the quantum isometry groups of these groups turn out to be 
  the free or tensor product of the quantum isometry groups of the factors $\Gamma_i$'s. Here is a brief list of groups for which we have computed the quantum isometry groups in this paper: 
  \begin{itemize}
  \item  All finitely generated abelian groups which do not have factors of the form $\IZ_2$ or $\IZ_4$  in the direct product decomposition into cyclic subgroups. 
  \item  Free product of all cyclic groups which do not have factors of the form $\IZ_2$ or $\IZ_4$. 
  \item Some special cases of direct or free product of cyclic groups having $\IZ_2$ or $\IZ_4$ as factors.
  \item Coxeter group.
  \item Baumslag-Solitar group.
  \item  Dihedral, Tetrahedral, Icosahedral and Octahedral groups.
  \end{itemize}

 \section{Quantum isometry group of $ C_r^*(\Gamma)$: existence and some generalities}
 We begin with a few basic definitions and facts about quantum isometry groups of spectral triples defined by Bhowmick and Goswami in \cite{qorient}. We denote the algebraic tensor product, spatial (minimal) $C^*$-tensor product and maximal $C^*$-tensor product by $\ot$, $\hat{\ot}$ and $\ot^{max}$  respectively. We'll use the leg-numbering notation. Let $\mathcal{H}$ be a complex Hilbert space, $\mathcal{K}(\mathcal{H})$ the $C^*$ algebra of compact operators on it, and $\mathcal{Q}$ a unital $C^*$ algebra. The multiplier algebra $\mathcal{M}(\mathcal{K}(\mathcal{H})\hat{\ot} \mathcal{Q})$ has two natural embeddings into $\mathcal{M}(\mathcal{K}(\mathcal{H})\hat{\ot} \mathcal{Q} \hat{\ot} \mathcal{Q})$, one obtained by extending the map $x \mapsto x \ot 1$ and the second one is obtained by composing this map with the flip on the last two factors. We will write $\omega^{12}$ and $\omega^{13}$ for the images of an element $\omega \in \mathcal{
M}(\mathcal{K}(\mathcal{H})\hat{\ot} \mathcal{Q})$ under these two maps respectively. We'll denote by $\mathcal{H} \bar{\ot} \mathcal{Q}$ the Hilbert $C^*$-module obtained by completing $\mathcal{H} \ot \mathcal{Q}$ with respect to the norm induced by the $\mathcal{Q}$ valued inner product $\langle\langle \xi \ot q, \xi^{\prime} \ot q^{\prime}\rangle\rangle \ := \langle\xi,\xi^{\prime}\rangle q^*q^{\prime}$, where $\xi,\xi^{\prime} \in \mathcal{H}$ and $q,q^{\prime} \in \mathcal{Q}$.\\

\subsection{Basic definitions}
\bdfn
A compact quantum group (CQG for short) is a pair $(\mathcal{Q},\Delta)$, where $\mathcal{Q}$ is a unital $C^*$-algebra and $\Delta : \mathcal{Q} \rightarrow \mathcal{Q} \hat{\ot} \mathcal{Q} $ is a  unital $C^*$-homomorphism satisfying the following two conditions:
\begin{enumerate} 
\item $(\Delta \ot id)\Delta = (id \ot \Delta)\Delta $ (co-associativity ).
\item Each of the linear spans of $\Delta(\mathcal{Q})(1 \ot \mathcal{Q})$ and that of $\Delta(\mathcal{Q})(\mathcal{Q} \ot 1)$ is norm-dense in $\mathcal{Q} \hat{\ot} \mathcal{Q}$.
\end{enumerate}
\edfn
A CQG morphism from $(\mathcal{Q}_1,\Delta_1)$ to another $(\mathcal{Q}_2,\Delta_2)$ is a unital $C^*$-homomorphism $\pi: \mathcal{Q}_1 \mapsto \mathcal{Q}_2$ such that $(\pi \ot \pi)\Delta_1=\Delta_2 \pi$.
\bdfn
$(\mathcal{Q}_1, \Delta_1)$ is called a quantum subgroup of $(\mathcal{Q}_2, \Delta_2)$ if there exists a surjective $C^*$-homomorphism $\eta$ from $\mathcal{Q}_2$ to $\mathcal{Q}_1$ such that $(\eta \ot \eta)\Delta_2=\Delta_1 \eta$ holds.   
\edfn
Sometimes we may denote the CQG $(\mathcal{Q},\Delta)$ simply as $\mathcal{Q}$, if $\Delta$ is clear from the context.
\bdfn
 A unitary (co)representation of $(\mathcal{Q},\Delta)$ on a Hilbert space $\mathcal{H}$ is a $\mathbb{C}$-linear map from $\mathcal{H}$ to the Hilbert module $\mathcal{H} \bar{\ot} \mathcal{Q}$ such that
\begin{enumerate} 
 \item $\langle\langle U(\xi),U(\eta)\rangle\rangle=<\xi,\eta>1_\mathcal{Q}$  (for all $\xi,\eta \in \mathcal{H}$). 
\item  $(U \ot id)U = (id \ot \Delta)U$.
\item Span $ \lbrace U(\xi)b: \ \xi \in \mathcal{H}, \ b \in \mathcal{Q}\rbrace $ is dense in $\mathcal{H} \bar{\ot} \mathcal{Q}$.  
\end{enumerate}
\edfn
Given such a unitary representation we have a unitary element $\tilde{U}$ belonging to $\mathcal{M}(\mathcal{K}(\mathcal{H})\hat{\ot} \mathcal{Q})$ given by $\tilde{U}(\xi \ot b)= U(\xi)b,(\xi \in \mathcal{H}, \ b \in \mathcal{Q})$ satisfying $(id \ot \Delta)(\tilde{U})=\tilde{U}^{12}\tilde{U}^{13}$.\\
\brmrk
It is known that the linear span of matrix elements of a finite dimensional unitary representation forms a dense Hopf *- algebra $\mathcal{Q}_0$ of $(\mathcal{Q},\Delta)$, on which an antipode $\kappa$ and co-unit $\epsilon$ are defined.
\ermrk
\bdfn
A closed subspace $\mathcal{H}_{1}$ of $\mathcal{H}$ is said to be invariant if $U(\mathcal{H}_{1}) \subseteq \mathcal{H}_{1} \bar{\ot} \mathcal{Q}$. A unitary representation $U$ of a CQG is said to be irreducible if there is no proper invariant subspace. 
\edfn
We denote by $\hat{\mathcal{Q}}$ the set of inequivalent irreducible representations of $(\mathcal{Q},\Delta)$. For $\pi \in \hat{\mathcal{Q}}$, let $d_{\pi}$ and $\lbrace t_{jk}^{\pi}: j, \ k = 1,\cdot\cdot\cdot,d_{\pi}\rbrace$ be the dimension and matrix coefficients of the corresponding finite dimensional representation respectively w.r.t some basis $\lbrace e_1,e_2,\cdot\cdot\cdot,e_{d_{\pi}}\rbrace$, i.e. $ \pi(e_i)= \sum_{j=1}^{d_{\pi}} e_j \ot t_{ij}^{\pi}$. Note that $\mathcal{Q}_0 =span\lbrace t_{ij}^{\pi}| \ \forall \ \pi \in \hat{\mathcal{Q}}  \rbrace$. The coproduct of $\mathcal{Q}$ is given by $\Delta(t_{ij}^{\pi})= \sum _{k=1}^{d_{\pi}} t_{kj}^{\pi} \ot t_{ik}^{\pi}$. Then for each $\pi \in \hat{\mathcal{Q}}$, we have a unique $   d_{\pi} \times d_{\pi}$ complex matrix $F_{\pi}$ such that
\begin{enumerate}
\item $F_{\pi}$ is positive and invertible with $Tr(F_{\pi})= Tr(F_{\pi}^{-1})= M_{\pi} > 0$ (say).
\item $h(t_{ij}^{\pi}t_{kl}^{\pi *}) = 1/M_{\pi} \ \delta_{ik}F_{\pi}(j,l)$.
\end{enumerate}
Corresponding to $\pi \in \hat{\mathcal{Q}}$, let $\rho_{sm}^{\pi}$ be the linear functional on $\mathcal{Q}$ given by $\rho_{sm}^{\pi}(x)= h(x_{sm}^{\pi}x), \ s,m= 1,\cdot\cdot\cdot,d_{\pi}$ for $x \in \mathcal{Q}$, where $x_{sm}^{\pi}= (M_{\pi})t_{km}^{\pi *}(F_{\pi})_{ks}$. Also let $\rho^{\pi}= \sum_{s=1}^{d_{\pi}} \rho_{ss}^{\pi}$.\\
\bdfn
 We say that a CQG $(\mathcal{Q},\Delta)$ acts on a unital $C^*$-algebra $B$ if there is a unital $C^*$-homomorphism (called action) $\alpha : B \rightarrow B \hat{\ot} \mathcal{Q} $ satisfying the following:
 \begin{enumerate}
 \item $(\alpha \ot id)\alpha = (id \ot \Delta)\alpha $.
 \item Linear span of $\alpha(B)(1 \ot \mathcal{Q})$ is norm-dense in $B \hat{\ot} \mathcal{Q}$.
 \end{enumerate}
\edfn
\bdfn
The action is said to be faithful if the $*$-algebra generated by the set $\lbrace (f \ot id)\alpha(b) \ \forall \ f \in B^{*}, \ \forall \ b\in B \rbrace$ is norm-dense in $\mathcal{Q}$, where $B^{*}$ is the Banach space dual of B.
\edfn
\brmrk
Given an action $\alpha$ of a CQG $\mathcal{Q}$ on a unital $C^*$-algebra B, we can always find a norm-dense, unital $*$-subalgebra $B_{0} \subseteq B$ such that $\alpha|_{B_{0}}: B_{0} \mapsto B_{0} \ot \mathcal{Q}_{0}$ is a Hopf-algebraic co-action. Moreover, $\alpha$ is faithful if and only if the $*$-algebra generated by $\lbrace (f \ot id)\alpha(b) \ \forall \ f \in B_{0}^{*}, \ \forall \ b\in B_{0} \rbrace$ is the whole of $\mathcal{Q}_{0}$.  
\ermrk
 Now we can define a projection $P_{\pi}: B \rightarrow B$ by $P_{\pi} := (id \ot \rho^{\pi})\alpha$ (note that $(id \ot \phi)\alpha(B)\subseteq B $ for all bounded linear functionals $\phi$ on $\mathcal{Q}$). We denote $Im (P_{\pi})$ by $B_{\pi}$, and call it the spectral subspace coming from $\pi$. We say that the action is of full spectrum if $B_{\pi} \neq 0 \ \forall \ \pi \in \hat{\mathcal{Q}}$.\\\\   
Given two CQG's $\mathcal{Q}_1$, $\mathcal{Q}_2$ the free product $\mathcal{Q}_1 \star \mathcal{Q}_2$   
as well as the maximal tensor product $\mathcal{Q}_1 \ot^{max} \mathcal{Q}_2$ admit the natural CQG structures, 
as given in \cite{free_wang}, \cite{tensor wang}. Moreover, they  have the following universal properties (see \cite{free_wang}, \cite{tensor wang}).   
\bppsn
\begin{enumerate}
\item  The canonical injections, say $i_1, \ i_2$ ($j_1, \ j_2$ respectively)  from  $\mathcal{Q}_1$ and $\mathcal{Q}_2$ to  $\mathcal{Q}_1 \star \mathcal{Q}_2$  
($\mathcal{Q}_1 \ot^{max} \mathcal{Q}_2$ 
respectively) are CQG morphisms. 
\item  Given  any CQG $\mathcal{C}$ and  morphisms $\pi_{1} : \mathcal{Q}_1 \mapsto \mathcal{C}$ and  
$\pi_{2} : \mathcal{Q}_2 \mapsto \mathcal{C}$  there always exists a unique morphism denoted by $\pi:=\pi_1 \ast \pi_2$ from  $\mathcal{Q}_1 \star \mathcal{Q}_2$ to  $\mathcal{C}$ satisfying
 $\pi \circ i_k=\pi_k$ for $k=1,2$. 
\item Furthermore, if the ranges of $\pi_{1}$ and $\pi_2$ commute, i.e.  $\pi_1(a)\pi_2(b)=\pi_2(b)\pi_1(a)$ $\forall \ a \in \mathcal{Q}_1, \ 
 b \in \mathcal{Q}_2$,  we  have a unique morphism $\pi^\prime$ from $\mathcal{Q}_1 \ot^{max} \mathcal{Q}_2$ to $\mathcal{C}$ satisfying $\pi^\prime \circ j_k=\pi_k$ for $k=1, 2.$
 \item The above conclusions hold for free or maximal tensor product of any finite number of CQG's as well.
 \end{enumerate}
\eppsn
\bdfn
Let $(\mathcal{A}^{\infty},\mathcal{H},\mathcal{D})$ be a spectral triple of compact type (a la Connes). Consider the category $Q(\mathcal{D})\equiv Q(\mathcal{A}^{\infty},\mathcal{H},\mathcal{D})$ whose objects are $(\mathcal{Q}, U)$, where $(\mathcal{Q},\Delta)$ is a CQG having a unitary representation U on the Hilbert space $\mathcal{H}$ satisfying the following:
\begin{enumerate}
\item $\tilde{U}$ commutes with $(\mathcal{D}\ot 1_{\mathcal{Q}})$.
\item $(id \ot \phi)\circ ad_{\tilde{U}}(a) \in (\mathcal{A}^{\infty})^{\prime\prime}$ for all $a \in \mathcal{A}^{\infty}$ and $\phi$ is any state on $\mathcal{Q}$, where $ad_{\tilde{U}}(x): = \tilde{U}(x \ot 1)\tilde{U}^*$ for $x \in \mathcal{B}(\mathcal{H})$. 
 \end{enumerate}
A morphism between two such objects  $(\mathcal{Q}, U)$ and $ (\mathcal{Q}^\prime, U^\prime)$ is a CQG morphism $\psi : \mathcal{Q} \rightarrow \mathcal{Q}^\prime $ such that $U^\prime = (id \ot \psi)U$. If a universal object exists in $Q(\mathcal{D})$ then we denote it by $\widetilde{QISO^{+}(\mathcal{A}^{\infty},\mathcal{H},\mathcal{D})}$ and the corresponding largest Woronowicz subalgebra for which $ad_{\tilde{U_0}}$ is faithful, where $U_{0}$ is the unitary representation of $\widetilde{QISO^{+}(\mathcal{A}^{\infty},\mathcal{H},\mathcal{D})}$, is called the quantum group of orientation preserving isometries and denoted by $QISO^{+}(\mathcal{A}^{\infty},\mathcal{H},\mathcal{D})$.
\edfn
Let us state Theorem $2.23$ of \cite{qorient} which gives a sufficient condition for the existence of  $QISO^{+}(\mathcal{A}^{\infty},\mathcal{H},\mathcal{D})$.
\bthm\label{existence thm}
 Let $(\mathcal{A}^{\infty},\mathcal{H},\mathcal{D})$ be a spectral triple of compact type. Assume that $\mathcal{D}$ has one dimensional kernel spanned by a vector $\xi \in \mathcal{H} $ which is cyclic and separating for $\mathcal{A}^{\infty}$ and each eigenvector of \ $\mathcal{D}$ belongs to $\mathcal{A}^{\infty}\xi$. Then QISO$^+(\mathcal{A}^{\infty},\mathcal{H},\mathcal{D})$ exists. 
 \ethm
 Let $(\mathcal{A}^{\infty},\mathcal{H},\mathcal{D})$ be a spectral triple satisfying the condition of Theorem \ref{existence thm} and $\mathcal{A}_{00}=Lin \lbrace a \in \mathcal{A}^{\infty} : a\xi$ is an eigenvector of $\mathcal{D} \rbrace$. Moreover, assume that $\mathcal{A}_{00}$ is norm-dense in $\mathcal{A}^{\infty}$. Let $\hat{\mathcal{D}}: \mathcal{A}_{00} \mapsto \mathcal{A}_{00}$ be defined by $\hat{\mathcal{D}}(a)\xi=\mathcal{D}(a\xi) (a \in \mathcal{A}_{00})$. This is well defined as $\xi$ is cyclic and separating vector for $\mathcal{A}^{\infty}$. Let $\tau$ be the vector state corresponding to the vector $\xi$.  
\bdfn\label{category definition}
Let $\mathcal{A}$ be a $C^*$-algebra and $\mathcal{A}^{\infty}$ be a dense *-subalgebra such that $(\mathcal{A}^{\infty},\mathcal{H},\mathcal{D})$ is a spectral triple as above. Let $\hat{\bold{C}}  (\mathcal{A}^{\infty},\mathcal{H},\mathcal{D})$ be the category with objects $(\mathcal{Q}, \alpha)$ such that $\mathcal{Q}$ is a CQG with a $C^*$-action $\alpha$ on $\mathcal{A}$ such that
\begin{enumerate}
\item $\alpha$ is $\tau$ preserving, i.e. $(\tau \ot id)\alpha(a)=\tau(a).1$ for all $a \in \mathcal{A}$.  
 \item $\alpha$ maps $\mathcal{A}_{00}$ into $\mathcal{A}_{00} \ot \mathcal{Q}$.
 \item $ \alpha\hat{\mathcal{D}}= (\hat{\mathcal{D}}\ot I)\alpha.$
\end{enumerate}
The morphisms in $\hat{\bold{C}}  (\mathcal{A}^{\infty},\mathcal{H},\mathcal{D})$ are CQG morphisms intertwining the respective actions.
 \edfn
 \bppsn{\label{new cat ppsn}}
 It is shown in Corollary 2.27 of \cite{qorient} that $QISO^{+}(\mathcal{A}^{\infty},\mathcal{H},\mathcal{D})$ is the universal object in  $\hat{\bold{C}}  (\mathcal{A}^{\infty},\mathcal{H},\mathcal{D})$.
 \eppsn 
 
 \subsection{QISO for a spectral triple on $C_r^*(\Gamma)$}
Now we discuss the special case of our interest. Connes considered this spectral triple in \cite{connes_compact}. 
Let $\Gamma$ be a finitely generated discrete group with generating set $S=\lbrace a_1,a_1^{-1}, a_2,a_2^{-1},$ $\cdot\cdot\cdot, a_k,a_k^{-1}\rbrace$. We make the convention of choosing the generating set to be symmetric, i.e. $a_i \in S$ implies  $a_i^{-1} \in S \ \forall \ i$. In case some $a_i$ has order 2, we include only $a_i$, i.e. not count it twice. 
 The corresponding word length function on the group defined by $l(g)=$ min $\lbrace r \in \mathbb{N}, \ g=h_1h_2\cdot\cdot\cdot h_r\rbrace$ where $h_i \in S$, i.e. for each i, $ h_i= a_j$ or $ a_j^{-1}$ for some $j$.
 Notice that $S=\lbrace g\in \Gamma| \ l(g)=1\rbrace$, using this length function we can define a metric on $\Gamma$ by $d(a,b)=l(a^{-1}b) \ \forall \ a,b \in \Gamma$. This is called the word metric corresponding to the generating set S.
Now consider the algebra $C_r^*(\Gamma)$, which is the $C^*$-completion of the group ring $\mathbb{C}\Gamma$ viewed as a subalgebra of $B(l^2(\Gamma))$ in the natural way via the left regular representation. We define a Dirac operator $D_{\Gamma}(\delta_g)=l(g)\delta_g$. In general, $D_\Gamma$ is an unbounded operator.
 $$Dom (D_\Gamma)=\lbrace \xi \in l^2(\Gamma): \sum_{g \in \Gamma} l(g)^2|\xi (g)|^2 < \infty \rbrace. $$
Here, $\delta_g$ is the vector in $l^{2}(\Gamma)$ which takes value $1$ at the point $g$ and $0$ at all other points. Natural generators of the algebra $\mathbb{C}\Gamma$ (images in the left regular representation ) will be denoted by $\lambda_g$, i.e. $\lambda_g (\delta_h)= \delta_{gh}$.
Let us define
$$\Gamma_r= \lbrace \delta_g| \ l(g)=r\rbrace,$$
$$\Gamma_{\leq r}= \lbrace \delta_g| \ l(g)\leq r\rbrace.$$
Moreover, let $p_r$ and $q_r$ be the orthogonal projections onto $Sp(\Gamma_r)$ and $Sp(\Gamma_{\leq r})$ respectively. Clearly 
$$D_{\Gamma}=\sum_{n \in \mathbb{N}_{0}} np_n,$$
where $p_r=q_r -q_{r-1}$ and $p_0=q_0$. The canonical trace on $C_r^*(\Gamma)$ is given by $\tau(\sum_{g\in \Gamma} c_g\lambda_g)=c_e$.  
It is easy to check that $ (\mathbb{C}\Gamma$, $l^2(\Gamma),D_\Gamma)$ is a spectral triple using Lemma 1.1 of \cite{oza_rie}.
Now take $\mathcal{A}= C_r^*(\Gamma), \ \mathcal{A}^\infty = \mathbb{C}\Gamma, \ \mathcal{H}= l^2(\Gamma)$ and $\mathcal{D}= D_{\Gamma}$ as before. Then QISO$^+(\mathbb{C}\Gamma$, $l^2(\Gamma),D_\Gamma)$ exists by Theorem \ref{existence thm}, taking $\delta_e$ as the cyclic separating vector for $ \mathbb{C}\Gamma$. As QISO$^+(\mathbb{C}\Gamma$, $l^2(\Gamma),D_\Gamma)$  depends on the generating set of $\Gamma$ it is denoted by $\mathbb{Q}(\Gamma,S)$. Most of the times we denote it by $\mathbb{Q}(\Gamma)$ if S is understood from the context.
Now as in \cite{grp algebra} its action $\alpha$ (say) on $C_r^*(\Gamma)$ is determined by  $$\alpha(\lambda_{\gamma})= \sum_{\gamma^\prime \in S}   \lambda_{\gamma^\prime} \ot q_{\gamma, \gamma^{\prime}} ,$$
where the matrix $[q_{\gamma,\gamma^\prime}]_{\gamma,\gamma^\prime \in S}$ is called the  fundamental representation in $M_{card(S)}(\mathbb{Q}$  $(\Gamma,S))$. Note that we have $\Delta(q_{\gamma,\gamma^\prime})= \sum_{\beta}q_{\beta, \gamma^{\prime}} \ot q_{\gamma, \beta}$. \\
$\mathbb{Q}(\Gamma,S)$ is also the universal object in the category $\hat{\bold{C}}(\mathbb{C}\Gamma$, $l^2(\Gamma),D_\Gamma)$ by Proposition \ref{new cat ppsn} and observe that all the eigenspaces of $\hat{\mathcal{D}_{\Gamma}}$, where $\hat{\mathcal{D}_{\Gamma}}$ is  as in Definition \ref{category definition}, are invariant under the action. The eigenspaces of $\hat{\mathcal{D}_{\Gamma}}$ are precisely $Span\lbrace \lambda_g| \ l(g)=r\rbrace$ with $r \geq 0$.   

It can also be identified with the universal object of some other categories naturally arising in the context. Consider the category $\bold{C}_{\bold{\tau}}$ of CQG's consisting of the objects $(\mathcal{Q},\alpha)$ such that $\alpha$ is an action of $\mathcal{Q}$ on $C_r^*(\Gamma)$ satisfying the following two properties:
\begin{enumerate}
\item $\alpha$ leaves $Sp(\Gamma_1)$ invariant.
 \item It preserves the canonical trace $\tau$ of $C_r^*(\Gamma)$.
 \end{enumerate}
Morphisms in $\bold{C}_{\bold{\tau}}$ are CQG morphisms intertwining the respective actions.
\blmma \label{category lemma}
The two categories $\bold{C}_{\bold{\tau}}$ and $\hat{\bold{C}}(\mathbb{C}\Gamma$, $l^2(\Gamma),D_\Gamma)$ are isomorphic.
\elmma 
For the proof the reader is referred to Lemma 2.16 of \cite{mandal 2}.\\
\bcrlre\label{corolary category}
It follows from Lemma \ref{category lemma} that there is a universal object, say $(\mathcal{Q}_{\tau}, \alpha_{\tau})$ in $\bold{C}_{\bold{\tau}}$ and $(\mathcal{Q}_{\tau},\alpha_{\tau}) \cong \mathbb{Q}(\Gamma,S)$.  
\ecrlre
We now identify $\mathbb{Q}(\Gamma,S)$ as a universal object in yet another category. 
Let us recall the quantum free unitary group $A_u(n)$ introduced in \cite{free_wang}. It is the universal unital $C^*$-algebra generated by $((a_{ij}))$ subject to the conditions that $((a_{ij}))$ and $((a_{ji}))$ are unitaries. Moreover, it admits a co-product structure with comultiplication $\Delta(a_{ij})=\Sigma_{l=1}^{n} a_{lj} \ot a_{il}$. Consider the category $\bold{C}$ with objects $(\mathcal{C},\lbrace x_{ij}, i,j=1,\cdot\cdot\cdot, 2k \rbrace )$ where $\mathcal{C}$ is a unital $C^*$-algebra generated by $((x_{ij}))$ such that $((x_{ij}))$ as well as $((x_{ji}))$ are unitaries and there is a unital $C^*$- homomorphism $\alpha_{\mathcal{C}}$  from $C_r^*(\Gamma)$ to $C_r^*(\Gamma) \hat{\ot} \ \mathcal{C}$ sending $e_i$ to $\sum_{j=1}^{2k} e_j \ot x_{ij}$, where $e_{2i-1}=\lambda_{a_{i}}$ and $e_{2i}=\lambda_{a_{i}}^{-1} \ \forall \ i=1,\cdot\cdot,k$. The morphisms from $(\mathcal{C},\lbrace x_{ij}, i,j=1,\cdot\cdot\cdot, 2k \rbrace )$ to $(\mathcal{P},\lbrace p_{ij}, i,j=1,\cdot\cdot\cdot, 2k \rbrace )$ are unital $*$-homomorphisms $\beta: \mathcal{C} \mapsto \mathcal{P}$ such that $\beta(x_{ij})=p_{ij}$.\\
Moreover, by definition of each object $(\mathcal{C},\lbrace x_{ij}, i,j=1,\cdot\cdot\cdot, 2k \rbrace )$ we get a unital $*$-morphism $\rho_{\mathcal{C}}$ from $A_u(2k)$ to $\mathcal{C}$ sending $a_{ij}$ to $x_{ij}$. Let the kernel of this map be $\mathcal{I}_{\mathcal{C}}$ and $\mathcal{I}$ be intersection of all such ideals. Then $\mathcal{C}^{\mathcal{U}}: = A_u(2k)/\mathcal{I}$ is the universal object generated by $x_{ij}^{\mathcal{U}}$ in the category $\bold{C}$. Furthermore, we can show, following a line of arguments similar to those in Theorem 4.8 of \cite{metric iso}, that it has a CQG structure with the co-product $\Delta(x_{ij}^{\mathcal{U}})=\sum_{l} x_{lj}^{\mathcal{U}}\ot x_{il}^{\mathcal{U}}$.      
\bppsn \label{proposition category}
$(\mathcal{Q}_{\tau},\alpha_{\tau})$ and $\mathcal{C}^{\mathcal{U}}$ are isomorphic as CQG.
\eppsn
{\it Proof:}\\
Let $((q_{ij}))$, $1\leq i,j \leq 2k$ be the fundamental representation of $(\mathcal{Q}_{\tau},\alpha_{\tau})$. Then we know that $((q_{ij}))$ and $((q_{ji}))$ are unitaries.  By the universal property of $\mathcal{C}^{\mathcal{U}}$, we always get a surjective map from $\mathcal{C}^{\mathcal{U}}$ to $(\mathcal{Q}_{\tau},\alpha_{\tau})$ sending $x_{ij}^{\mathcal{U}}$ to $q_{ij}$, which intertwins the actions too. \\
On the other hand, we can construct a state on $C_r^*(\Gamma)$ defined by $\tilde{\tau}: = (\tau \ot h)\alpha_{\mathcal{C}}^{\mathcal{U}}$, where h is the Haar state of $(\mathcal{C}^{\mathcal{U}},x_{ij}^{\mathcal{U}})$, is clearly $\alpha_{\mathcal{C}}^{\mathcal{U}}$ invariant. Note that $(C_r^*(\Gamma),\Delta_{\Gamma}, t_{ij})$ is an object in the category $\bold{C}$, where $((t_{ij}))$ is the diagonal matrix with entries $t_{2i-1,2i-1}=\lambda_{a_i}, \ t_{2i,2i}=\lambda_{{a_i}^{-1}}$. By the universal property of $\mathcal{C}^{\mathcal{U}}$ we always get a surjective *- morphism $\pi$ from $\mathcal{C}^{\mathcal{U}}$ to $C_r^*(\Gamma)$ sending            $x_{ij}^{\mathcal{U}}$ to $t_{ij}$ such that $ (id \ot \pi)\alpha_{\mathcal{C}}^{\mathcal{U}}=\Delta_{\Gamma}$ holds. Moreover, we get 
$$(\tilde{\tau} \ot id)\Delta_{\Gamma}(x)=(id \ot \tilde{\tau})\Delta_{\Gamma}(x)=\tilde{\tau}(x).1  \ \forall \ x \in C_r^*(\Gamma). $$
But we know that the canonical trace (Haar state) $\tau$ is the unique bi-invariant state on $(C_r^*(\Gamma), \Delta_{\Gamma})$. Hence, $\tilde{\tau}=\tau$. Then by universality of $(\mathcal{Q}_{\tau},\alpha_{\tau})$ we get a surjective map from $(\mathcal{Q}_{\tau},\alpha_{\tau})$ to   $\mathcal{C}^{\mathcal{U}}$ sending $q_{ij}$ to $x_{ij}^{\mathcal{U}}$, which also intertwins the actions. Thus, $(\mathcal{Q}_{\tau},\alpha_{\tau})$ and $\mathcal{C}^{\mathcal{U}}$ are isomorphic as CQG.    
\qed\\
If $\Gamma$ is commutative, the maximal commutative quantum subgroup of $\mathbb{Q}(\Gamma)$ (its abelianization) is denoted by $C(ISO(\Gamma))$. In Section $3$ we will see many examples for which $\mathbb{Q}(\Gamma)\cong C(ISO(\Gamma))$.
Now we fix some notational conventions which will be useful in later sections.
Note that the action $\alpha$ is of the form 
\begin{eqnarray*}
 \alpha(\lambda_{a_{1}}) &=& \lambda_{a_{1}} \ot A_{11} + \lambda_{a_{1}^{-1}} \ot A_{12} + \lambda_{a_{2}} \ot A_{13} + \lambda_{a_{2}^{-1}} \ot A_{14} + \cdot \cdot \cdot     + \\ 
 &&\lambda_{a_{k}} \ot A_{1(2k-1)} + \lambda _{a_{k}^{-1}} \ot A_{1(2k)}, \\
  \alpha(\lambda_{a_{1}^{-1}}) &=& \lambda_{a_{1}} \ot A_{12}^* + \lambda_{a_{1}^{-1}} \ot A_{11}^* + \lambda_{a_{2}} \ot A_{14}^* + \lambda_{a_{2}^{-1}} \ot A_{13}^* + \cdot \cdot \cdot     + \\ 
  &&\lambda_{a_{k}} \ot A_{1(2k)}^* + \lambda _{a_{k}^{-1}} \ot A_{1(2k-1)}^*, \\
 \alpha(\lambda_{a_{2}}) &=& \lambda_{a_{1}} \ot A_{21} + \lambda_{a_{1}^{-1}} \ot A_{22} + \lambda_{a_{2}} \ot A_{23} + \lambda_{a_{2}^{-1}} \ot A_{24} + \cdot \cdot \cdot     + \\
&& \lambda_{a_{k}} \ot A_{2(2k-1)} + \lambda _{a_{k}^{-1}} \ot A_{2(2k)}, \\
\alpha(\lambda_{a_{2}^{-1}}) &=& \lambda_{a_{1}} \ot A_{22}^* + \lambda_{a_{1}^{-1}} \ot A_{21}^* + \lambda_{a_{2}} \ot A_{24}^* + \lambda_{a_{2}^{-1}} \ot A_{23}^* + \cdot \cdot \cdot     + \\ 
&&\lambda_{a_{k}} \ot A_{2(2k)}^* + \lambda _{a_{k}^{-1}} \ot A_{2(2k-1)}^*, \\
\vdots  &&   \hspace{1cm}           \vdots \\
\alpha(\lambda_{a_{k}}) &=& \lambda_{a_{1}} \ot A_{k1} + \lambda_{a_{1}^{-1}} \ot A_{k2} + \lambda_{a_{2}} \ot A_{k3} + \lambda_{a_{2}^{-1}} \ot A_{k4} + \cdot \cdot \cdot     + \\ 
 &&\lambda_{a_{k}} \ot A_{k(2k-1)} + \lambda _{a_{k}^{-1}} \ot A_{k(2k)}, \\
 \alpha(\lambda_{a_{k}^{-1}}) &=& \lambda_{a_{k}} \ot A_{k2}^* + \lambda_{a_{1}^{-1}} \ot A_{k1}^* + \lambda_{a_{2}} \ot A_{k4}^* + \lambda_{a_{2}^{-1}} \ot A_{k3}^* + \cdot \cdot \cdot     + \\ 
  &&\lambda_{a_{k}} \ot A_{k(2k)}^* + \lambda _{a_{k}^{-1}} \ot A_{k(2k-1)}^*. \\ 
\end{eqnarray*} 
 From this we get the unitary corepresentation \\
 \begin{equation}\label{eq matrix}
 U \equiv ((u_{ij}))= 
 \begin{pmatrix}
 A_{11} & A_{12} & A_{13} & A_{14} & \cdots & A_{1(2k-1)} & A_{1(2k)}\\
 A_{12}^* & A_{11}^* & A_{14}^* & A_{13}^* & \cdots &  A_{1(2k)}^* & A_{1(2k-1)}^*\\
 A_{21} & A_{22} & A_{23} & A_{24} & \cdots & A_{2(2k-1)} & A_{2(2k)}\\
 A_{22}^* & A_{21}^* & A_{24}^* & A_{23}^* & \cdots & A_{2(2k)}^* & A_{2(2k-1)}^*\\
 \vdots  && \hspace{1cm}        \vdots \\
 A_{k1} & A_{k2} & A_{k3} & A_{k4} & \cdots & A_{k(2k-1)} & A_{k(2k)}\\
 A_{k2}^* & A_{k1}^* & A_{k4}^* & A_{k3}^* & \cdots &  A_{k(2k)}^* & A_{k(2k-1)}^*\\ 
 \end{pmatrix}.
 \end{equation}
 
From now on, we call it as fundamental unitary.
The coefficients $A_{ij}$ and $A_{ij}^*$'s generate a norm-dense subalgebra of $\mathbb{Q}(\Gamma,S)$. We also note that the antipode of $\mathbb{Q}(\Gamma,S)$ maps $u_{ij}$ to $u_{ji}^*$.
\begin{rmrk}\label{another descrip}
Using Corollary \ref{corolary category} and Proposition \ref{proposition category},  $\mathbb{Q}(\Gamma,S)$ is the universal unital $C^*$-algebra generated by $A_{ij}$ as above subject to the relations that U, $U^t$ are unitaries and $\alpha$ given above is a $C^*$-homomorphism on $C_r^*(\Gamma)$.
\end{rmrk}
 \bppsn\label{prop 1}
$(C^*(\Gamma),\Delta_{\Gamma})$ always acts on $C^*_r(\Gamma)$ isometrically and faithfully, i.e. $(C^*(\Gamma),\Delta_{\Gamma})$ is a subobject of $\mathbb{Q}(\Gamma,S)$ in the category $\hat{\bold{C}}(\mathbb{C}\Gamma$, $l^2(\Gamma),D_\Gamma)$.  
\eppsn 
{\it Proof:}\\
The usual (co)action of $ C^*(\Gamma)$ on $C^*_r(\Gamma)$ coming from the coproduct of $C^*(\Gamma)$ gives us action on itself given by 
$\Delta_{\Gamma}(\lambda_g)=\lambda_g \ot \lambda_g, \forall \ g\in \Gamma$. It is clear that this action  is isometric and faithful. 
So by Remark \ref{another descrip}, there is a surjective $C^*$-morphism from $\mathbb{Q}(\Gamma)$ to $C^*(\Gamma)$ sending the entries $A_{i(2i-1)}, \ A_{i(2i-1)}^*$ of the fundamental unitary mentioned before to $\lambda_{a_i}$ and $\lambda _{a_{i}^{-1}}$ respectively, others being sent to zero. The morphism intertwins the actions too.         \qed \\\\
The above proposition tells us that for a nonabelian group $\Gamma$, the quantum isometry group $\mathbb{Q}(\Gamma)$ is always a genuine CQG, i.e. the underlying $C^*$-algebra is non-commutative.
\bcrlre \label{group algebra cor}
$\mathbb{Q}(\Gamma,S) \cong (C^*(\Gamma),\Delta_{\Gamma}) $ if and only if the matrix (\ref{eq matrix}) is diagonal.  
\ecrlre
{\it proof:}\\
Suppose $\mathbb{Q}(\Gamma,S) \cong (C^*(\Gamma),\Delta_{\Gamma})$, then clearly (\ref{eq matrix}) is a diagonal matrix. Conversely, let the matrix (\ref{eq matrix}) is diagonal. By Proposition \ref{prop 1} we always get a surjective $C^*$-morhism from $\mathbb{Q}(\Gamma,S)$ to $(C^*(\Gamma),\Delta_{\Gamma})$ sending $A_{i(2i-1)}$ to $\lambda_{a_{i}}$ for each i. Now the action of $\mathbb{Q}(\Gamma,S)$ is defined as $\alpha(\lambda_{a_{i}})= \lambda_{a_{i}} \ot A_{i(2i-1)}$. So we get for $g=a_{i_1}a_{i_2}\cdots a_{i_k}, \ \alpha(\lambda_g)=\lambda_g \ot A_g$ by the $*$-homomorphism property of $\alpha$, where $A_g \in \mathbb{Q}(\Gamma,S)$ is defined by $A_g := A_{i_1(2i_1-1)} \cdots A_{i_k(2i_k-1)}$. It is easy to see that the map $g \mapsto A_g$ satisfies $A_{gh}=A_{g}A_{h}, \ A_{g^{-1}}=A_{g}^*, \ A_e=1$. Then by the universal property of $(C^*(\Gamma),\Delta_{\Gamma})$ we get a surjective $*$-morphism from $(C^*(\Gamma),\Delta_{\Gamma}) $ to $\mathbb{Q}(\Gamma,S)$ sending $\lambda_{a_{i}}$ to $A_{i(2i-1)}$. This completes the proof.          
\qed\\\\
We end the discussion of this subsection with the following easy observation which will be used in various places of this article. 
\bppsn\label{normal comp}
If $UV=0$ for two normal elements in a $C^*$-algebra then 
$$ U^*V=VU^*=0,$$
$$ V^*U=UV^*=VU=0.$$ 
\eppsn

\subsection{Recollection of some known facts}
To the best of our knowledge, first computations of  $\mathbb{Q}(\Gamma)$ were done in \cite{grp  algebra}. Thereafter, several articles by different authors were 
 devoted to computations of $\mathbb{Q}(\Gamma)$ for concrete groups.  In \cite{grp algebra} quantum isometry groups of cyclic groups (except $\mathbb{Z}_4$) 
 were shown to be commutative.  In case of  $\mathbb{Z}_4$ it turns out to be  noncommutative and infinite dimensional.
It is in fact isomorphic with $C^*(D_{\infty} \times \mathbb{Z}_2)$ as a $C^*$-algebra (see  \cite{free cyclic}).
Later in \cite{org filt} it was identified with $\mathbb{Z}_2 \wr_{\ast} \mathbb{Z}_2$ as a quantum group.
 
 The authors of \cite{S_n} introduced the doubling procedure, and moreover they showed that for 
 the symmetric group $S_n$ with standard generating sets consisting of $(n-1)$ transpositions, the quantum isometry group coincides with 
  the doubling of the group algebra. The same result holds for $D_{2(2n+1)}$ as well (see \cite{dihedral}). 
  We will also briefly discuss  the  doubling procedure in Subsection \ref{doub section}.

 In \cite{two par}, Banica and Skalski introduced  two parameter families $H^+_s(p,q)$, 
 of quantum symmetry groups and they studied quantum isometry groups of duals of free product of cyclic groups for several cases in \cite{free cyclic}. They showed that
 $$H_n^{+}\cong \mathbb{Q}(\underbrace{\mathbb{Z}_2 \ast \mathbb{Z}_2 \cdot\cdot\cdot\ast  \mathbb{Z}_2)}_{n \ copies},$$
$$ H^{+}(n,0)\cong \mathbb{Q}(\underbrace{\mathbb{Z} \ast \mathbb{Z} \cdot\cdot\cdot\ast  \mathbb{Z})}_{n \ copies},$$
 $$ H^{+}_s(n,0)\cong \mathbb{Q}(\underbrace{\mathbb{Z}_s \ast \mathbb{Z}_s \cdot\cdot\cdot\ast  \mathbb{Z}_s)}_{n \ copies},$$ where $s\neq 2,4$.
It was observed in \cite{org filt}  $$C(O_n^{-1})\cong \mathbb{Q}(\underbrace{\mathbb{Z}_2 \times \mathbb{Z}_2 \cdot\cdot\cdot\times  \mathbb{Z}_2)}_{n \ copies}.$$

\subsection{The case when $\Gamma$ is a free or direct product } 
\bthm\label{prop 4}
Let $\Gamma_1,\Gamma_2,\cdot\cdot\cdot, \Gamma_k$ be finitely generated discrete groups with the symmetric generating sets $S_1,S_2,\cdot\cdot\cdot,S_k$ respectively. Consider $\Gamma=\Gamma_1 \times \Gamma_2 \times \cdot\cdot\cdot \times \Gamma_k$ with the generating set $S=\cup_{i=1}^k S^{\prime}_i$, where $S^{\prime}_i=(0,\cdot\cdot,S_i,\cdot\cdot,0) \ \forall \ i=1,\cdot\cdot,k$. Then $\mathbb{Q}(\Gamma,S)$ has $\mathbb{Q}(\Gamma_1,S_1)\ot^{max} \mathbb{Q}(\Gamma_2,S_2)\ot^{max}\cdot\cdot\cdot \ot^{max} \mathbb{Q}(\Gamma_k,S_k)$ as a quantum subgroup.
\ethm
{\it Proof:}\\
Let $S_i=\lbrace a_{ij}, \ j=1,2,\cdot\cdot k_i\rbrace$ where $i=1,2,\cdot\cdot,k$ and $((u_{pj}^{(i)}))$ be the fundamental unitary representation for the action of $\mathbb{Q}_i \equiv \mathbb{Q}(\Gamma_i,S_i)$ on $C_r^*(\Gamma_i)$. Consider the action $\alpha$ :  $C_r^*(\Gamma)\mapsto  C_r^*(\Gamma) \hat{\ot} \ \mathbb{Q}$ given by
 $$\alpha(\lambda_{a_{ip}^{\prime}})= \Sigma_{j=1}^{k_i} \lambda_{a_{ij}^{\prime}} \ot 1_{(1)} \ot 1_{(2)}\cdot\cdot\ot 1_{(i-1)} \ot u_{pj}^{(i)} \ot 1_{(i+1)}\ot \cdot\cdot \ot 1_{(k)},$$  where $\mathbb{Q}=\mathbb{Q}_1 \ot^{max} \mathbb{Q}_2 \ot^{max} \cdot\cdot \ot^{max} \mathbb{Q}_k $ and $1_{(l)}$ denotes the identity element of the underlying $C^*$-algebra of $\mathbb{Q}_l$. It is easy to verify that this gives an isometric action of $\mathbb{Q}$ on $C_r^*(\Gamma)$, hence by the universality of $\mathbb{Q}(\Gamma,S)$ we get a surjective morphism from $\mathbb{Q}(\Gamma,S)$ to $\mathbb{Q}$. \qed \\   
\brmrk\label{prop 5}
In the set up of the previous theorem, replace the direct product by free product, i.e. take $\Gamma=\Gamma_1 \ast \Gamma_2 \ast \cdot\cdot\cdot \ast \Gamma_k$. Then $\mathbb{Q}(\Gamma_1) \star \mathbb{Q}(\Gamma_2) \star \cdot\cdot\cdot \star \mathbb{Q}(\Gamma_k)$ is a quantum subgroup of $\mathbb{Q}(\Gamma)$. The proof is very similar to Theorem \ref{prop 4} and hence omitted. 
\ermrk
\brmrk
We give an example to show that $\mathbb{Q}(H)$ may not be a quantum subgroup of $\mathbb{Q}(\Gamma)$ in general for a subgroup H of $\Gamma$, when $\Gamma$ is neither $H \times K$ nor $H \ast K$ for some $K$. 
Take $S_4$ with generating sets $(12),(23),(34)$ and $H$ be the subgroup of it defined by $ <(12),(34)> $. $H$ is clearly isomorphic with $\mathbb{Z}_2 \times \mathbb{Z}_2 $. We know that $\mathbb{Q}(\mathbb{Z}_2 \times \mathbb{Z}_2)$ is infinite dimensional (see \cite{org filt}) but the underlying $C^*$-algebra of $\mathbb{Q}(S_4)$ is isomorphic to $C^*(S_4)\oplus C^*(S_4) $ (result in \cite{S_n}). 
\ermrk
\brmrk
We'll be usually interested to see whether $\mathbb{Q}(\Gamma_1) \star \mathbb{Q}(\Gamma_2) \star \cdot\cdot\cdot \star \mathbb{Q}(\Gamma_k)$ or  $\mathbb{Q}(\Gamma_1) \ot^{max} \mathbb{Q}(\Gamma_2) \ot^{max} \cdot\cdot\cdot \ot^{max} \mathbb{Q}(\Gamma_k)$ coincides with $\mathbb{Q}(\Gamma)$ (in section $3$ and $4$) whenever $\Gamma$ is either $\Gamma_1 \ast \Gamma_2 \ast \cdot\cdot\cdot \ast \Gamma_k$ or $\Gamma_1 \times \Gamma_2 \times \cdot\cdot\cdot \times \Gamma_k$ respectively. In this context the following observation will be useful.
\ermrk
\blmma\label{lemma in sec 1}
Let $((u_{ij}))$ be the fundamental representation of $\mathbb{Q}(\Gamma)$ on $C_r^*(\Gamma)$. Also assume $\Gamma=\Gamma_1 \ast \Gamma_2 \ast \cdot\cdot\cdot \ast \Gamma_k$. Then $\mathbb{Q}(\Gamma)\cong \mathbb{Q}(\Gamma_1) \star \mathbb{Q}(\Gamma_2) \star \cdot\cdot\cdot \star \mathbb{Q}(\Gamma_k)$ if $((u_{ij}))$ is of the block diagonal form 
$$U=\begin{pmatrix}
C_1 & 0 \cdots & 0\\
0 & C_2 \cdots & 0\\
\vdots & \vdots & \vdots\\
0 & 0 \cdots & C_k\\
\end{pmatrix}$$
with respect to the decomposition of the generating set S of $\Gamma$ into $S_1\cup S_2\cup\cdot\cdot\cup S_k$.   
\elmma
{\it Proof:}\\
Let us write $S_i=\lbrace a_{ij}, \ j=1,2,\cdot\cdot,k_i\rbrace$ as in Theorem \ref{prop 4}. It is clear from the above form of U that for each i, the action (say $\alpha$ ) of $\mathbb{Q}(\Gamma)$ maps $C_r^*(\Gamma_i)\cong C^*(\lambda_{a_{ij}}, \ j=1,2,\cdot\cdot,k_i)$ to $C_r^*(\Gamma_i) \hat{\ot} \mathbb{Q}(\Gamma,S)$ with the corresponding fundamental unitary being $C_i=((c_{pj}^{(i)}))$. This means we have a morphism from $\mathbb{Q}_i\equiv\mathbb{Q}(\Gamma_i,S_i)$ to $\mathbb{Q}(\Gamma,S)$ sending $u_{pj}^{(i)}$ to $c_{pj}^{(i)}$. By definition of the free product, this gives a morphism from $\mathbb{Q}_1 \star \mathbb{Q}_2 \star \cdot\cdot \star \ \mathbb{Q}_k$ to $\mathbb{Q}(\Gamma)$. By Remark \ref{prop 5} we always get a surjective morphism from $\mathbb{Q}(\Gamma)$ to $\mathbb{Q}_1 \star \mathbb{Q}_2 \star \cdot\cdot \star \ \mathbb{Q}_k$ sending $c_{pj}^{(i)}$ to $u_{pj}^{(i)} \ \forall \ i$. This completes the proof. \qed \\ 
\brmrk\label{rmk in sec 1}
The above lemma is true if we replace the free product by direct product, i.e.  $\mathbb{Q}(\Gamma)\cong \mathbb{Q}(\Gamma_1) \ot^{max} \mathbb{Q}(\Gamma_2) \ot^{max} \cdot\cdot\cdot \ot^{max} \mathbb{Q}(\Gamma_k)$ if and only if $((u_{ij}))$ is of the block diagonal form. Moreover, entries of one such block commute with the entries of any other block. The proof of this fact is very similar to Lemma \ref{lemma in sec 1}, hence omitted.  
\ermrk
\subsection{$\mathbb{Q}(\Gamma)$ as a doubling of $C^*(\Gamma)$}\label{doub section}
We briefly recall the doubling procedure of the group algebra from \cite{S_n}, \cite{doubling}.
Let $(\mathcal{Q},\Delta)$ be a CQG with a CQG-automorphism $\theta$ such that $\theta^2=id$. The doubling of this CQG, 
 say $(\mathcal{D}_{\theta}(\mathcal{Q}),\tilde{\Delta})$, is  given by  $\mathcal{D}_{\theta}(\mathcal{Q}) :=\mathcal{Q}\oplus \mathcal{Q}$ (direct sum as a $C^*$-algebra),
 and the coproduct is defined by the following, where we have denoted  the injections of $\mathcal{Q}$ onto 
 the first and second coordinate in $\mathcal{D}_{\theta}(\mathcal{Q})$ by $\xi$ and $\eta$ respectively, i.e. 
$\xi(a)=(a,0), \  \eta(a)= (0,a), \ (a \in \mathcal{Q}).$
$$\tilde{\Delta} \circ \xi= (\xi \ot \xi + \eta \ot [\eta \circ \theta])\circ  \Delta,$$
$$\tilde{\Delta} \circ \eta= (\xi \ot \eta + \eta \ot [\xi \circ \theta])\circ \Delta.$$
Below we give a sufficient condition for the quantum isometry group to be a doubling of the group algebra. For this, it is convenient to use a slightly different notational convention:
 let $U_{2i-1,j}=A_{ij}$ for $i=1,\ldots, k,$ \ $j=1, \ldots, 2k$ and $U_{2i,2l}= A_{i(2l-1)}^*, \  U_{2i,2l-1}= A_{i(2l)}^*$ for  $i=1,\ldots, k,$ \ $l=1, \ldots, k$.  
\blmma \label{doub lemma}
Let $\Gamma$ be a group with $k$ generators $\lbrace a_1,a_2,\cdot\cdot\cdot, a_k\rbrace$ and define $\gamma_{2l-1}:=a_l,~ \gamma_{2l}:=a^{-1}_l \ \forall \ l=1,2,\cdot\cdot\cdot,k$. Now $\sigma$ be an order 2 automorphism  on the set  $\lbrace 1,2,\cdot\cdot,2k-1,2k\rbrace$ and $\theta$ be an automorphism of the group given by $\theta(\gamma_i)=\gamma_{\sigma(i)}$. We assume the following:
\begin{enumerate}
 \item  $B_i: = U_{i,\sigma(i)}\neq 0 \  \forall \ i$ and  $ U_{i,j}=0 \  \forall \ j \not\in \lbrace \sigma(i),i\rbrace$,
\item  $A_iB_j=B_jA_i=0 \ \forall \ i,j$  such that $\sigma(i)\neq i$ and  $\sigma(j)\neq j,$   where    $A_i=U_{i,i}$,
 \item  All   $U_{i,j}U_{i,j}^*$   are   central  projections,
 \item  There   are   well  defined  $C^*$-isomorphisms   $\pi_1,\pi_2$  from  $C^*(\Gamma)$  to $ C^*\lbrace A_i, i=1,2,\cdot\cdot,2k \rbrace$   and   $C^*\lbrace B_i, i=1,2,\cdot\cdot,2k \rbrace$  respectively  such  that  
 $$\pi_1(\lambda_{a_i})=A_i, \ \pi_2(\lambda_{a_i})=B_i \ \forall \ i.$$
\end{enumerate}
 Then $\mathbb{Q}(\Gamma)$ is doubling of $(C^*(\Gamma,\Delta_{\Gamma})$ corresponding to a given automorphism $\theta$. Moreover, the fundamental unitary takes the following form 
$$ \begin{pmatrix}
 A_{1} & 0 & 0 & 0 & \cdots & 0 & B_{1}\\
 0 & A_{2} & 0 & 0 & \cdots &  B_{2} & 0\\
 0 & 0 & A_{3} & 0 & \cdots & 0 & 0\\
 0 & 0 & 0 & A_{4} & \cdots & 0 & 0\\
 \vdots  && \hspace{1cm}        \vdots \\
 0 & B_{2k-1} & 0 & 0 & \cdots & A_{2k-1} & 0\\
 B_{2k} & 0 & 0 & 0 & \cdots &  0 & A_{2k}\\ 
 \end{pmatrix}.$$\\
\elmma
{\it Proof:}\\
First note that when $\sigma(i)=i$, then i-th row contains only one non-zero element $A_i$, then $U_{i,i}: = A_i=B_i$. Without loss of generality we assume that there are first $k_1$ number of i's from the set $\lbrace 1,2,\cdot\cdot,2k-1,2k\rbrace$ such that $\sigma(i) \neq i$ and $k_2$ number of j's from the set $\lbrace k_1+1,\cdot\cdot,2k\rbrace$ such that $\sigma(j) = j$. This implies $k_1 + k_2=2k$ and always $k_1 \geq 1$ as $\sigma$ is non-trivial.  
  Now the $C^*$-algebra  $C^*(\Gamma) \oplus C^*(\Gamma)$ is generated by $(\lambda_{\gamma_i} \oplus 0), (0 \oplus \lambda_{\gamma_{\sigma(i)}})$, $i=1,2, \ldots,2k$.   
   We can define a $C^*$-homomorphism $\pi$ from $\mathcal{D}_{\theta}(C^*(\Gamma))$ to $\mathbb{Q}(\Gamma)$ given by 
   $$\pi(\lambda_{\gamma_i} \oplus 0)= A_i,~~~ \pi(0 \oplus \lambda_{\gamma_{\sigma(i)}})= B_i \ \forall \ i=1,2,\cdot\cdot k_1,$$ 
  $$ \pi(\lambda_{\gamma_j} \oplus 0)= A_jA_1A_1^{*}, ~~~ \pi(0 \oplus \lambda_{\gamma_{\sigma(j)}})= A_jB_1B_1^{*} \ \forall \ j=k_1+1,\cdot\cdot 2k.$$
  It is easy to verify that this is indeed a CQG isomorphism.
 \qed
\subsection{$\mathbb{Q}(\Gamma,S)$ as a quantum isometry group of metric space}
In this subsection, our aim is to identify $\mathbb{Q}(\Gamma,S)$ with the quantum isometry group of some metric space in the sense of \cite{metric iso}, in case $\Gamma$ is abelian. 
We first recall the definition of quantum isometry group in the purely metric space setting. Given a compact metric space $(X,d)$, we say that an action $\alpha$ 
of a CQG on $C(X)$ is isometric if the action $\alpha_r = (id \ot \pi_r)\circ \alpha$ of the reduced CQG $\mathbb{Q}_r$ 
(where $\pi_r : \mathbb{Q} \rightarrow \mathbb{Q}_r$ is the canonical map from $\mathbb{Q}$ to the reduced CQG  $\mathbb{Q}_r$) satisfies $$\alpha_r(d_x)(y)= \kappa (\alpha_r(d_y)(x)),$$
$\forall \ x,y \in X$, where $d_x(.)\equiv d(x,.)$ and $\kappa$ denotes the (norm bounded by Theorem 3.23 in \cite{huang}) antipode of $\mathbb{Q}_r$. It is shown in \cite{metric iso} that 
in case $X \subseteq \mathbb{R}^n$ isometrically embedded, the above condition is equivalent to the following: 
$$\sum_i (F_i(x)-F_i(y))^2 = \sum_i (x_i-y_i)^2.1,$$ 
where $F_i(x)= \alpha(X_i)(x)$ and  $X_i$ denotes the $i$-th coordinate function of $\mathbb{R}^n$, restricted to X. 
It can be easily seen, by almost verbatim adaptation of the arguments in \cite{metric iso}, that a similar result would hold if we replaced
$\mathbb{R}^n$ by $\mathbb{C}^n$. That is, for $X\subseteq \mathbb{C}^n$ isometrically, with the metric $d(z,w)^2 = \Sigma_i|z_i-w_i|^2$, 
 a CQG action $\alpha : C(X)\rightarrow C(X) \hat{\ot} \ \mathbb{Q}$ is isometric in the metric space sense if and only if $$\sum_i(F_i(z)-F_i(w))^*(F_i(z)-F_i(w))= d^2(z,w)1,$$
$\forall \ z,w \in X$. Moreover, it is clear that a sufficient condition for the above is that 
$$ \sum_i F_i(z)^*F_i(w)= \langle z,w \rangle 1 \equiv \sum_i \bar{z_i}w_i1.$$
We can also prove, as in \cite{metric iso}, that for metric spaces $(X,d)$ isometrically embedded in $\mathbb{C}^n$, there exists a universal CQG acting isometrically on it,
to be denoted by QISO$(X,d)$, and its action is affine in the sense that  $\alpha(Z_i)= \Sigma_j Z_j \ot b_{ij} + 1 \ot R_i $ for some $b_{ij}$ and $R_i$, where $Z_i$ is the coordinate function restricted to $X \subseteq \mathbb{C}^n$.\\  
Let us now consider a finitely generated abelian group $\Gamma$  with a  symmetric generating set $S=\lbrace\gamma_1,\gamma_2\cdot\cdot\cdot,\gamma_n\rbrace$ 
and let $((q_{ij}))$ be the fundamental unitary for $\mathbb{Q}(\Gamma,S)$. 
Now consider the dual group of $\Gamma$ say $G=\hat{\Gamma}$, which is a compact topological group.
Moreover, $\hat{\Gamma}$ can be identified with a compact subset X of $\mathbb{C}^n$ via the map $\chi \mapsto (\chi(\gamma_1),\chi(\gamma_2),\cdot\cdot\cdot,\chi(\gamma_n))$. 
There is a natural Euclidean metric $\hat{d_S}$ on $\hat{\Gamma}$ given by $\hat{d_S}(\chi,\chi ^\prime)^2= \Sigma_i |\chi(\gamma_i) -\chi^\prime(\gamma_i)|^2$.
\bthm
With the above set-up, we have  $$\mathbb{Q}(\Gamma,S)\cong QISO(\hat{\Gamma},\hat{d_S}).$$
  
\ethm
{\it Proof:}\\
It is well-known that $C_r^*(\Gamma)$ is isomorphic with $C(\hat{\Gamma})$ via the Fourier transform $\mathcal{F}$,  extended as a unitary from $\ell^2(\Gamma)$ to $L^2(\hat{\Gamma})$ i.e. $C(\hat{\Gamma})=\mathcal{F}C_r^*(\Gamma)\mathcal{F}^{-1}\subseteq B(L^2(\Gamma))$. 
Let $U$ denote the unitary representation of $\mathbb{Q}(\Gamma,S)$ on $C_r^*(\Gamma)$ and let 
 $U^\prime$ and $\alpha^\prime$ be the corresponding representation and action of $\mathbb{Q}(\Gamma,S)$ on $C(\hat{\Gamma})$ respectively, i.e.
 $U^\prime=\mathcal{F}U\mathcal{F}^{-1}$ and $\alpha^{\prime}(Z_{i})=(\mathcal{F}\ot 1)\alpha(Z_i)(\mathcal{F}^{-1}\ot 1).$    
As $\lambda_{\gamma_1},\lambda_{\gamma_2},\cdot\cdot,\lambda_{\gamma_n}$ generate  $C_r^*(\Gamma)$ as a $C^*$-algebra, 
under the Fourier transform $C(\hat{\Gamma})$ is generated by  $\hat{\lambda}_{\gamma_1},\hat{\lambda}_{\gamma_2},$ $\cdot\cdot,\hat{\lambda}_{\gamma_n}$.
We have $Z_i=\hat{\lambda}_{\gamma_i}$ as  the coordinate functions restricted to $\hat{\Gamma}\subset\mathbb{C}^n$. 
Clearly, $C(\hat{\Gamma})\cong C^*\lbrace Z_1,Z_2,\cdot\cdot,Z_n \rbrace \subseteq C_0(\mathbb{C}^n)$.
Let $F(z)=(F_1(z),F_2(z),$ $\cdot\cdot,F_n(z))$ for $z \in \hat{\Gamma}$, where $F_i(z):= \alpha^\prime(Z_i)(z)$.  
Now it is clear that $F_i\equiv \alpha^\prime(Z_i)=\Sigma_j Z_j \ot q_{ij}$,
and as $((q_{ij}))$ is unitary, it is easy to verify that  $ \Sigma_i F_i(z)^*F_i(w)= \Sigma_i \bar{z_i}w_i \ot 1_{\mathbb{Q}}$.
Thus $\alpha^\prime$ is an isometric action of $\mathbb{Q}(\Gamma,S)$ on $C(\hat{\Gamma})$, i.e. $\mathbb{Q}(\Gamma,S)$ is a quantum subgroup of QISO$(\hat{\Gamma},\hat{d_S})$ using the surjective morphism $b_{ij} \mapsto q_{ij}, R_i \mapsto 0$.

Conversely, assume that the isometric action of $\mathbb{Q}^\prime\equiv$ QISO$(\hat{\Gamma},\hat{d_S})$ is given by $\beta(Z_i)= \Sigma_j Z_j \ot b_{ij} + 1 \ot R_i $ 
for some $b_{ij}$ and $R_i$. 
Consider any Borel, probability measure on $(\hat{\Gamma},\hat{d_S})$ and convolving with the Haar state of $\mathbb{Q}^\prime$, we get
a $\mathbb{Q}^\prime$ invariant probability measure say $\mu$ on $\hat{\Gamma}$. As $\hat{\Gamma}$ itself acts isometrically on $(\hat{\Gamma},\hat{d_S})$, 
$C_r^*(\Gamma)(\cong C(\hat{\Gamma})$) is a quantum subgroup of $\mathbb{Q}^\prime$, hence $\mu$ is a $\hat{\Gamma}$ invariant probability measure on $\hat{\Gamma}$.  
Therefore $\mu$ must be the (unique) Haar measure of $\hat{\Gamma}$. 
As the state corresponding to the Haar measure (say $\phi$) maps each $Z_i$ to zero,
we conclude using the $\beta$ invariance of $\phi$ i.e. $(\phi \ot id)\beta(Z_i)= \phi(Z_i)1_{\mathbb{Q}^\prime}$, that $R_i=0$ for each $i$. 
Thus $\beta$ is linear, i.e. $\beta(Z_i)=\Sigma_j Z_j \ot b_{ij}$. 
The corresponding $*$-homomorphism from $C_r^*(\Gamma)$ to $C_r^*(\Gamma) \hat{\ot} \mathbb{Q}^\prime$ maps 
$\lambda_{\gamma_i}$ to $\Sigma_j \lambda_{\gamma_j}\ot b_{ij}$ and it is clearly an action. 
 Finally, as $\beta$ preserves the Haar state $\mu$, it extends to a unitary representation on the $L^2$-space, hence 
 in particular, $((b_{ij}))$ is the matrix corresponding to a  unitary representation.  In other words, $(( b_{ij}))$ is a unitary element of $M_n({\mathbb{Q}}^\prime).$
  This shows by the definition of $\mathbb{Q}(\Gamma,S)$ that $\mathbb{Q}^\prime$ is a quantum subgroup of $\mathbb{Q}(\Gamma,S)$, 
  with the surjective morphism sending $q_{ij}$ to $b_{ij}$. This completes the proof. \ \qed\\

  \subsection{Polynomial growth of $\mathbb{Q}(\Gamma)$}
  We briefly discuss some sufficient conditions for the quantum group $\mathbb{Q}(\Gamma)$ to have polynomial growth property in the sense of \cite{jyotish_new}, when $\Gamma$
   has polynomial growth.
   

        We now state and prove the main result of this subsection.
        \bthm
      Let $\Gamma$ be a finitely generated discrete group with polynomial growth and assume that the action of $\mathbb{Q}(\Gamma)$ on $C^*_r(\Gamma)$ has full spectrum. 
      Then the dual discrete quantum group of $\mathbb{Q}(\Gamma)$ also has polynomial growth property.
        \ethm
        {\it Proof:}\\
       Let $S$ be a finite generating set for $\Gamma$. As the group $\Gamma$ has polynomial growth, there is some polynomial $p$ of one variable such that 
        the cardinality of the set $\{ g_1g_2\ldots g_m:~g_i \in S \ \forall \ i, ~m \leq n \}$ is bounded 
        by $p(n)$ for each $n$. That is, the dimension of the  vector space, say $\clv_n$, 
        spanned by elements of the form $\lambda_{a_1} \ldots \lambda_{a_m}$ where $a_i \in S$ and $m \leq n$, has dimension less than or equal to $p(n)$. But this space is clearly left invariant 
         by the action of $\mathbb{Q}(\Gamma)$. Let us denote the restriction of this action to $\clv_n$, which is a finite dimensional unitary representation, by 
          $\pi_n$. Moreover, by the assumption of full spectrum, every irreducible representation of $\mathbb{Q}(\Gamma)$ must be a sub-representation of some $\pi_n$ for sufficiently 
           large $n$. This allows us to define a central length function $l$ on the set of irreducible representation of $\mathbb{Q}(\Gamma)$ (in the sense of \cite{jyotish_new}) 
             by setting $l(\pi)$ to be equal to the smallest value of $n$ for which $\pi$ is a sub-representation of $\pi_n$. Clearly,
             $$ \sum_{ \pi :~l(\pi) \leq n } {d_\pi}^2 \leq {p(n)},$$ where $d_\pi$ denotes the dimension of the irreducible $\pi$. From this, it is easily seen that this length function 
              satisfies the criteria of Definition 4.1 of \cite{jyotish_new}, hence the dual of $\mathbb{Q}$  has polynomial growth. \qed\\

\subsection{The structure of the maximal commutative subgroup of $\mathbb{Q}(\Gamma)$ for  $\Gamma = \mathbb{Z}_{n}^k$}  
\bppsn \label{iso prop}
Let $\Gamma=\underbrace{(\mathbb{Z}_n \times \mathbb{Z}_n \times\cdot\cdot \mathbb{Z}_n)}_{k \ copies}$ and $S_k$ be the group of permutation of k elements. We have:\\
$1.$ If $n=2$ then we have $C(ISO(\Gamma))\cong$ $ C(\widehat{\underbrace{(\mathbb{Z}_2 \times \mathbb{Z}_2 \times\cdot\cdot \mathbb{Z}_2)}_{k \ copies}}$ $\rtimes S_k)$.\\
$2.$ If $n\neq 2, 4$ then $C(ISO(\Gamma)) \cong C(\widehat{\underbrace{(\mathbb{Z}_n \times \mathbb{Z}_n \times\cdot\cdot \mathbb{Z}_n)}_{k \ copies}}$ $\rtimes (\mathbb{Z}_2^k \rtimes S_k)). $
\eppsn
{\it Proof:}\\
Let $a_1,a_2,\cdot\cdot\cdot,a_k$ be the usual generating elements for $\Gamma$. 
The fundamental unitary is of the form\\ 
$ U \equiv ((u_{ij}))=
 \begin{pmatrix}
 A_{11} & A_{12} & A_{13} & A_{14} & \cdots & A_{1(2k-1)} & A_{1(2k)}\\
 A_{12}^* & A_{11}^* & A_{14}^* & A_{13}^* & \cdots &  A_{1(2k)}^* & A_{1(2k-1)}^*\\
 A_{21} & A_{22} & A_{23} & A_{24} & \cdots & A_{2(2k-1)} & A_{2(2k)}\\
 A_{22}^* & A_{21}^* & A_{24}^* & A_{23}^* & \cdots & A_{2(2k)}^* & A_{2(2k-1)}^*\\
 \vdots  && \hspace{1cm}        \vdots \\
 A_{k1} & A_{k2} & A_{k3} & A_{k4} & \cdots & A_{k(2k-1)} & A_{k(2k)}\\
 A_{k2}^* & A_{k1}^* & A_{k4}^* & A_{k3}^* & \cdots &  A_{k(2k)}^* & A_{k(2k-1)}^*\\ 
 \end{pmatrix}$.\\\\ However, we have the additional condition that $u_{ij}$'s commute among themselves and each of them is normal. Our first claim is  $A_{ij}A_{ik}=0=A_{ji}A_{ki}, \  \forall \  i,j,k$ with $j \neq k$. \\\\
 We break it into two cases.
\paragraph{Case 1: $n=2$:}
Consider the term  $\alpha(\lambda_{a_i^2})= \lambda_e \ot 1_{\mathbb{Q}} \ \forall \ i$. Now comparing the coefficients of $\lambda_{a_la_m} \ \forall \  l\neq m$ on both sides of the equation we obtain $A_{il}A_{im}=0 \  \forall \  l\neq m$. Applying the antipode we find $A_{li}A_{mi}=0 \  \forall \  l\neq m$. 
 \paragraph{Case 2: $n\neq 2, 4$:} 
 Using the relation  $ \alpha(\lambda_{a_{j}}) \alpha(\lambda_{a_{j}^{-1}})=\alpha(\lambda_{a_{j}^{-1}}) \alpha(\lambda_{a_{j}})=\alpha(\lambda_e)= \lambda_e \ot 1_{\mathbb{Q}} $ and comparing the coefficients of 
$ \lambda_{a_{i}^{2}}$ and $ \lambda_{a_{i}^{-2}}$ on both sides we must have
\begin{equation} \label{eq 2.1}
  A_{i(2j-1)}A_{i(2j)}^*=A_{i(2j)}^*A_{i(2j-1)}=0.
 \end{equation} 
   Applying the antipode on (\ref{eq 2.1}) we get that 
 $A_{i(2j-1)}^*A_{i(2j)}^*=A_{i(2j)}^*A_{i(2j-1)}^*  =0$, which shows
 \begin{equation} \label{eq 2.2}
   A_{i(2j)}A_{i(2j-1)}=A_{i(2j-1)}A_{i(2j)}=0 \ \forall \ i,j.
  \end{equation} 
    Now, we will show that $A_{i(2j)}A_{im}=A_{i(2j-1)}A_{im}=0$, where $m \neq 2j,(2j-1)$. Using the condition $ \alpha(\lambda_{a_{i}}) \alpha(\lambda_{a_{i}^{-1}})=\alpha(\lambda_{a_{i}^{-1}}) \alpha(\lambda_{a_{i}})=\alpha(\lambda_e)= \lambda_e \ot 1_{\mathbb{Q}} $ and comparing the coefficients of $\lambda_{a_ja_l}, \ \lambda_{a_ja_l^{-1}}$ where $j\neq l$ one can get
    \begin{equation} \label{eq 2.3}
A_{i(2j-1)}A_{i(2l)}^* + A_{i(2l-1)}A_{i(2j)}^* =0,
\end{equation}
\begin{equation} \label{eq 2.4}
 A_{i(2j)}A_{i(2l-1)}^* + A_{i(2l)}A_{i(2j-1)}^*=0.
 \end{equation}
Multiplying $A_{i(2j-1)}^*$ and $A_{i(2j)}^*$ on the right side of the equations (\ref{eq 2.3}) and (\ref{eq 2.4}) respectively we have   
  $A_{i(2j-1)}A_{i(2l)}^*A_{i(2j-1)}^*=A_{i(2j)}A_{i(2l-1)}^*A_{i(2j)}^*= 0$. Thus, $A_{i(2j-1)}A_{i(2l)}^*=A_{i(2j)}A_{i(2l-1)}^*=0$ by using the $C^*$-norm condition and commutativity among $A_{ij}$'s. Now using the Proposition \ref{normal comp} we get $A_{i(2j-1)}A_{i(2l)}=A_{i(2j)}A_{i(2l-1)}=0$. Similarly, we have $A_{i(2j)}A_{i(2l)}=A_{i(2j-1)}A_{i(2l-1)}$ $ =0$. This proves the claim by using antipode. Moreover, for finite $n$, we have  $A_{ij}^*=A_{ij}^{n-1}$. \\\\
Now we produce the explicit isomorphism. Identify $C(\widehat{\underbrace{(\mathbb{Z}_n \times \mathbb{Z}_n \times\cdot\cdot \mathbb{Z}_n)}_{k \ copies}}\rtimes (\mathbb{Z}_2^k \rtimes S_k))$ as a $C^*$-algebra in a natural way with $C\widehat{\underbrace{(\mathbb{Z}_n \times \mathbb{Z}_n \times\cdot\cdot \mathbb{Z}_n)}_{k \ copies}} \ot  C(\mathbb{Z}_2^k) \ot  C(S_k)$. 
Let $\chi_{ij} \in C(S_k)$ be the characteristic function of the set of those permutations which map i to j and also $\eta_{0,i}, \ \eta_{1,i} \in C(\mathbb{Z}_2^k)$ be the characteristic functions of sets which have respectively 0 or 1 in the i-th coordinate. One can easily check as in Theorem $6.1$ of \cite{two par} that the map
$$A_{i(2j-1)} \mapsto z_i \ot \chi_{ij} \ot \eta_{0,i},$$
$$A_{i(2j)} \mapsto z_i \ot \chi_{ij} \ot \eta_{1,i},$$
gives an isomorphism between $C(ISO(\Gamma))$ and $C(\widehat{\underbrace{(\mathbb{Z}_n \times \mathbb{Z}_n \times\cdot\cdot \mathbb{Z}_n)}_{k \ copies}}\rtimes (\mathbb{Z}_2^k \rtimes S_k))$, and that it preserves the respective coproducts. So it becomes a CQG isomorphism.  \qed  \\ 

 \section{Computations for free product of cyclic groups}\label{free sec}
 \bthm\label{thm 1 for free}
 Let $\Gamma=\Gamma_1\ast\Gamma_2\cdot\cdot\cdot\ast \ \Gamma_l$ where $\Gamma_i=\underbrace{(\mathbb{Z}_{n_i}\ast\mathbb{Z}_{n_i}\cdot\cdot\ast\mathbb{Z}_{n_i})}_{k_i \ copies}$. Also assume   $n_1\neq n_2\neq\cdot\cdot \neq n_l$ and $n_i\neq2,4 \ \forall \ i$, then $\mathbb{Q}(\Gamma)$ will be $H^{+}_{n_{1}}(k_1,0)\star H^{+}_{n_{2}}(k_2,0)\star \cdot\cdot\star H^{+}_{n_{l}}(k_l,0) $.  i.e.  $\mathbb{Q}(\Gamma) \cong \mathbb{Q}(\Gamma_1)\star \mathbb{Q}(\Gamma_2)\star \cdot\cdot\cdot\star \mathbb{Q}(\Gamma_l)$. 
 \ethm
 {\it Proof:}\\
For simplicity of notation we present the case when all $k_i=1$. The general case will follow by essentially the same arguments, with just on extra bit of careful book-keeping of notations. So let $\Gamma= \mathbb{Z}_{n_1}\ast \mathbb{Z}_{n_2}\ast\cdot\cdot\cdot \ast \mathbb{Z}_{n_l}$ and let $a_i$ be the standard generators of $\mathbb{Z}_{n_i}, \ o(a_i)=n_i \ \forall \ i$. Take $\lbrace a_1, a_1^{-1}, a_2, a_2^{-1}\cdot\cdot\cdot, a_l, a_l^{-1} \rbrace$ as the symmetric generating set for this group. Without loss of generality we assume $n_1 < n_2 < \cdot\cdot\cdot < n_l$. $n_l$ can take the value $\infty$ also. Now the fundamental unitary of $\mathbb{Q}(\Gamma)$ is of the form \\
 $$ U=
 \begin{pmatrix}
 A_{11} & A_{12} & A_{13} & A_{14} & \cdots & A_{1(2l-1)} & A_{1(2l)}\\
 A_{12}^* & A_{11}^* & A_{14}^* & A_{13}^* & \cdots &  A_{1(2l)}^* & A_{1(2l-1)}^*\\
 A_{21} & A_{22} & A_{23} & A_{24} & \cdots & A_{2(2l-1)} & A_{2(2l)}\\
 A_{22}^* & A_{21}^* & A_{24}^* & A_{23}^* & \cdots & A_{2(2l)}^* & A_{2(2l-1)}^*\\
 \vdots  && \hspace{1cm}        \vdots \\
 A_{l1} & A_{l2} & A_{l3} & A_{l4} & \cdots & A_{l(2l-1)} & A_{l(2l)}\\
 A_{l2}^* & A_{l1}^* & A_{l4}^* & A_{l3}^* & \cdots &  A_{l(2l)}^* & A_{l(2l-1)}^*\\ 
 \end{pmatrix}.$$\\\\
 
 Our aim is to show that it reduces to the form\\\\
 
$$ \begin{pmatrix}
 A_{11} & A_{12} & 0 & 0 & \cdots & 0 & 0\\
 A_{12}^* & A_{11}^* & 0 & 0 & \cdots &  0 & 0\\
 0 & 0 & A_{23} & A_{24} & \cdots & 0 & 0\\
 0 & 0 & A_{24}^* & A_{23}^* & \cdots & 0 & 0\\
 \vdots  && \hspace{1cm}        \vdots \\
 0 & 0 & 0 & 0 & \cdots & A_{l(2l-1)} & A_{l(2l)}\\
 0 & 0 & 0 & 0 & \cdots &  A_{l(2l)}^* & A_{l(2l-1)}^*\\ 
 \end{pmatrix}$$\\\\
i.e. only the diagonal $(2\times2)$ block will survive and others become zero. By Lemma \ref{lemma in sec 1} this will complete the proof.\\
 We break the proof into a number of lemmas.\\
 \begin{lmma}\label{lemma for free}
 $A_{i(2j)}A_{i(2j-1)}= A_{i(2j-1)}A_{i(2j)}=0 \ \forall \ i,j.$
 \end{lmma}
 {\it Proof:}\\
  Using the relation  $ \alpha(\lambda_{a_{j}}) \alpha(\lambda_{a_{j}^{-1}})=\alpha(\lambda_{a_{j}^{-1}}) \alpha(\lambda_{a_{j}})=\alpha(\lambda_e)= \lambda_e \ot 1_{\mathbb{Q}} $ and comparing the coefficients of 
$ \lambda_{a_{i}^{2}}$ and $ \lambda_{a_{i}^{-2}}$ on both sides we have
 $A_{i(2j-1)}A_{i(2j)}^*=A_{i(2j)}^*A_{i(2j-1)}=0$. Applying the antipode one can get that 
 $A_{i(2j-1)}^*A_{i(2j)}^*=A_{i(2j)}^*A_{i(2j-1)}^*$ $=0 $, thus $ A_{i(2j)}A_{i(2j-1)}=A_{i(2j-1)}$ $A_{i(2j)}=0$.  \qed \\
  \begin{lmma}\label{lemma for notations}
 $A_{1j}=0 \  \forall \  j>  2. $
 \end{lmma}
 {\it Proof:}\\
First we fix some notational convention for clarity of exposition. 
Let $b_{2i-1}=a_i$ and $ b_{2i}=a_{i}^{-1}$, $i=1,2,\cdot\cdot,l$. Now consider $ \alpha(\lambda_{a_{1}^{n_{1}-1}})= \alpha(\lambda_{a_{1}^{-1}})$. Observe that both sides must contain only the words of length 1, i.e. only $b_{i}$'s occur in the chain. For $ \underline{i}= (i_1,i_2,\cdot\cdot\cdot,i_m)$, set $A_{\underline{i}}\equiv A_{1i_{1}}A_{1i_{2}}\cdot\cdot A_{1i_{m}}, \ b_{\underline{i}}\equiv b_{i_{1}}b_{i_{2}}\cdot\cdot b_{i_{m}} $ and $\Lambda =\lbrace  (i_1,i_2,\cdot\cdot\cdot,i_{n_1-1})| \  i_j \in \lbrace 1,2,\cdot\cdot,2l\rbrace\rbrace$, note that L.H.S $=\sum_{\underline{i}\in \Lambda} \lambda_{b_{\underline{i}}} \ot A_{\underline{i}}$.
 Using Lemma \ref{lemma for free} we conclude that $A_{\underline{i}}=0$ whenever $b_{\underline{i}}$ is not a reduced word, i.e. there is some $i_r$ such that $b_{i_{r}}^{-1}= b_{i_{r+1}}$. Thus the L.H.S reduces to $\sum_{\underline{i}\in \Lambda} \lambda_{b_{\underline{i}}} \ot A_{\underline{i}}$, where $b_{\underline{i}}$ is a reduced word. Our claim is that there are no coefficients of $ \lambda_{a_{d}}$ and $ \lambda_{a_{d}^{-1}} \ \forall \ d > 1$ in L.H.S. Suppose $a_{j_{1}}^{m_1}a_{j_{2}}^{m_2}\cdot\cdot a_{j_{p}}^{m_p}= a_2$ where $j_{i} \neq j_{i+1} \ \forall \ i$ and $ |m_1|+|m_2|+\cdot\cdot+|m_p|= n_1 -1$  which implies $a_{j_{1}}^{m_1}a_{j_{2}}^{m_2}\cdot\cdot\cdot a_{j_{p}}^{m_p}a_{2}^{-1}=e$. This shows that  p must be 1 and $j_1=2$. Hence we can write $a_{2}^{m_1}=a_2, \ |m_1|= n_1-1$, leading to a contradiction as $n_1 < n_2$. We get similar contradiction if we replace $a_2$ by any other $a_i \ \forall \ i > 2$ as $n_1 < n_i \ \forall \ i > 2$. Thus comparing the coefficients of $ \lambda_{a_{d}}$ and $ \lambda_{a_{d}^{-1}} \ \forall \
 d > 1$ on both sides one can deduce $A_{1j}^* =0 \ \forall \ j> 2$. Taking the adjoint we get the desired result. \qed \\  
 Now applying the antipode we find that  $A_{i1}=A_{i2}=0 \ \forall \ i> 1$.\\\\
 The structure of the unitary matrix reduces to the form \\
$$ \begin{pmatrix}
 A_{11} & A_{12} & 0 & 0 & \cdots & 0 & 0\\
 A_{12}^* & A_{11}^* & 0 & 0 & \cdots &  0 & 0\\
 0 & 0 & A_{23} & A_{24} & \cdots & A_{2(2l-1)} & A_{2(2l)}\\
 0 & 0 & A_{24}^* & A_{23}^* & \cdots & A_{2(2l)}^* & A_{2(2l-1)}^*\\
 \vdots  && \hspace{1cm}        \vdots \\
 0 & 0 & A_{l3} & A_{l4} & \cdots & A_{l(2l-1)} & A_{l(2l)}\\
 0 & 0 & A_{l4}^* & A_{l3}^* & \cdots &  A_{l(2l)}^* & A_{l(2l-1)}^*\\ 
 \end{pmatrix}.$$\\\\
 If we repeat this line of reasoning with other generators starting with $a_{2}$ and so on, we get the desired block diagonal form and hence complete the proof. One can easily observe that the last condition  $ \alpha(\lambda_{a_{l}^{n_{l}-1}})= \alpha(\lambda_{a_{l}^{-1}})$ is not used in the proof. So we can include the $n_l=\infty$ case also.   \qed\\\\
 We'll now show how to extend Theorem \ref{thm 1 for free} in some cases when $2$ or $4$ can occur in {$n_1,n_2,\cdot\cdot\cdot,n_l$}.
 \bthm \label{thm 2 for free}
 Theorem \ref{thm 1 for free} is valid also for the case when $n_1=2$ and $2\neq n_2\neq\cdot\cdot \neq n_l$ where $n_i\neq 4,\infty \ \forall \ i$.
  \ethm
  {\it Proof:}\\
  As in the proof of Theorem \ref{thm 1 for free} we assume $k_i=1 \ \forall \ i$ for simplicity of exposition and moreover $n_{2} < n_{3} \cdot\cdot\cdot < n_{l}$.
From the relation $\alpha(\lambda_{a_{1}^2})=\lambda_e \ot 1_{\mathbb{Q}}$, comparing the coefficients of $\lambda_{a_{2}^{2}},\lambda_{a_{2}^{-2}},\lambda_{a_{3}^{2}},\lambda_{a_{3}^{-2}},\cdot\cdot \lambda_{a_{l}^{2}},\lambda_{a_{l}^{-2}}$ on both sides we find that
\begin{equation} \label{eq 3.4.1}  
A_{12}^2=A_{13}^2= \cdot\cdot\cdot =A_{1(2l-2)}^2=A_{1(2l-1)}^2=0.
\end{equation} 
 This implies
 \begin{equation}\label{eq 3.4.2}
   A_{21}^2=( A_{21}^*)^2= \cdot\cdot\cdot =  A_{k1}^2=( A_{k1}^*)^2=0
  \end{equation}
  by applying the antipode on (\ref{eq 3.4.1}). Also observe that 
  \begin{equation}\label{eq 3.4.3}
  A_{2(2j+1)}A_{2(2j)}= A_{2(2j)}A_{2(2j+1)}=0 \ \forall \ j
  \end{equation}
  by the similar argument of Lemma \ref{lemma for free}.\\ 
   Now with the notations of Lemma \ref{lemma for notations}, consider $ \alpha(\lambda_{a_{2}^{n_{2}-1}})= \alpha(\lambda_{a_{2}^{-1}})$. 
L.H.S contains terms of the form $\sum_{\underline{i}\in \Lambda} \lambda_{b_{\underline{i}}} \ot A_{\underline{i}}$ 
where $\Lambda =\lbrace  (i_1,i_2,\cdot\cdot\cdot,i_{n_2-1})  | \  i_j \in \lbrace 1,2,\cdot\cdot,(2l-1)\rbrace\rbrace$. 
Note that here also  $b_{\underline{i}}$ is reduced word by the equations (\ref{eq 3.4.2}), (\ref{eq 3.4.3}) and the fact
$n_{2} < n_{3} \cdot\cdot\cdot < n_{l}$.
 Thus comparing the coefficients of $\lambda_{a_{1}},\lambda_{a_{p}}, \lambda_{a_{p}^{-1}} \ \forall \ p>2$ on both sides of the relation  $ \alpha(\lambda_{a_{2}^{n_{2}-1}})= \alpha(\lambda_{a_{2}^{-1}})$ one can deduce that \ 
$A_{21}^*=A_{2k}^*=0 \ \forall \ k>3$. \  
Applying the antipode we have, $A_{12}=A_{13}=0$ and $A_{k2}=A_{k3}=0 \ \forall \ k > 2$. \\
So the fundamental Unitary  is reduced to the form \\
$$ \begin{pmatrix}
 A_{11} & 0 & 0  & \cdots & A_{1(2l-2)} & A_{1(2l-1)} \\
 0 & A_{22} & A_{23}  & \cdots & 0 & 0\\
 0 & A_{23}^* & A_{22}^*  & \cdots & 0 & 0\\
 \vdots  && \hspace{1cm}        \vdots \\
 A_{l1} & 0 & 0 &  \cdots & A_{l(2l-2)} & A_{l(2l-1)}\\
 A_{l1}^* & 0 & 0  & \cdots &  A_{l(2l-1)}^* & A_{l(2l-2)}^*\\ 
 \end{pmatrix}.$$\\\\
Repeating similar arguments with $a_3,a_4,\cdot\cdot, a_l$ we will get the desired block diagonal form and complete the proof. \ \qed \\
\brmrk
We see that finiteness of order of $a_l$ is used in the last paragraph of Theorem \ref{thm 2 for free}. But in Theorem  $2.1$ this condition is not used. In fact it won't be necessary for any of the results that follow. So this is the only case where the finiteness condition is necessary.
\ermrk
\bthm\label{thm 3 for free}
The result of Theorem \ref{thm 1 for free} remains true for $n_1=4, \ 4\neq n_2\neq\cdot\cdot \neq n_l$ and $n_i\neq2 \ \forall \ i$.   
\ethm
  {\it Proof:}\\
  We continue using the notation and convention of Theorem \ref{thm 1 for free} and without loss of generality assume $k_i=1 \ \forall \ i$ and $o(a_1)=4$. First consider the term $ \alpha(\lambda_{a_{1}^{2}})$ and note that the coefficient of $\lambda_e$ in the expression of  $ \alpha(\lambda_{a_{1}^{2}})$ must be zero.\\
Thus we have
\begin{equation}\label{eq 3.6.1}  
  A_{12}A_{11}+A_{11}A_{12} + \cdot\cdot\cdot +  A_{1(2l-1)}A_{1(2l)}+A_{1(2l)}A_{1(2l-1)} =0.
 \end{equation}
 Now, observe that
 \begin{equation} \label{eq 3.6.2}
 A_{1(2i-1)}A_{1(2i)}= A_{1(2i)}A_{1(2i-1)}= 0 \ \forall \ i > 1
 \end{equation}
  by the same argument of Lemma \ref{lemma for free}. From (\ref{eq 3.6.1}) and (\ref{eq 3.6.2}) we have
  \begin{equation}\label{eq 3.6.3} 
  A_{12}A_{11}+A_{11}A_{12}=0. 
 \end{equation}
 Next, we want to show that $A_{1i}=0 \ \forall \ i >2$. From the equation $ \alpha(\lambda_{a_{1}^{3}})= \alpha(\lambda_{a_{1}^{-1}})$ comparing the coefficients of $\lambda_{a_{2}}, \lambda_{a_{2}^{-1}},
\lambda_{a_{3}}, \lambda_{a_{3}^{-1}},\cdot\cdot\cdot, \lambda_{a_{l}}, \lambda_{a_{l}^{-1}}$ on both sides we deduce \\
$$A_{1(2i)}^*=A_{1(2i-1)}( A_{12}A_{11}+A_{11}A_{12} + \cdot\cdot + A_{1(2i-1)}A_{1(2i)}+\cdot\cdot +A_{1(2l)}A_{1(2l-1)}),$$
$$A_{1(2i-1)}^*=A_{1(2i)}( A_{12}A_{11}+A_{11}A_{12}+ \cdot\cdot+ A_{1(2i)}A_{1(2i-1)}+\cdot\cdot +A_{1(2l)}A_{1(2l-1)}),$$ where $ \ \forall \ i> 1$. \\
Now, using (\ref{eq 3.6.2}) and (\ref{eq 3.6.3}) one can conclude that   $A_{1i}=0 \ \forall \ i >2$. The fundamental unitary is of the form \\
$$ \begin{pmatrix}
 A_{11} & A_{12} & 0 & 0 & \cdots & 0 & 0\\
 A_{12}^* & A_{11}^* & 0 & 0 & \cdots &  0 & 0\\
 0 & 0 & A_{23} & A_{24} & \cdots & A_{2(2l-1)} & A_{2(2l)}\\
 0 & 0 & A_{24}^* & A_{23}^* & \cdots & A_{2(2l)}^* & A_{2(2l-1)}^*\\
 \vdots  && \hspace{1cm}        \vdots \\
 0 & 0 & A_{l3} & A_{l4} & \cdots & A_{l(2l-1)} & A_{l(2l)}\\
 0 & 0 & A_{l4}^* & A_{l3}^* & \cdots &  A_{l(2l)}^* & A_{l(2l-1)}^*\\ 
 \end{pmatrix}.$$\\\\
 We then follow the strategy of Theorem \ref{thm 1 for free} to complete the proof. \qed\\
\bthm\label{thm 4 for free}
The conclusion of Theorem \ref{thm 1 for free} holds also for $n_1=2, n_2=4$ where  $ n_3\neq n_4\neq\cdot\cdot \neq n_l$ and $n_i\neq 2, 4, \infty \ \forall \ i\geq 3$.
\ethm
  {\it Proof:}\\
With the notation and convention as before, we get from the condition $ \alpha(\lambda_{a_{2}^{3}})=  \alpha(\lambda_{a_{2}^{-1}})$ comparing the coefficient of $ \lambda_{a_{1}}$ on both sides
\begin{equation} \label{eq 3.7.1}
A_{21}^*=A_{21}(A_{22}A_{23}+A_{23}A_{22})
\end{equation}
as $A_{2(2i+1)}A_{2(2i)}=A_{2(2i)}A_{2(2i+1)}=0 \ \forall \ i > 1 $  by  arguments similar to those of Lemma \ref{lemma for free}.\\
Again, considering the coefficient of $\lambda_e$ in the expression of $ \alpha(\lambda_{a_{2}^{2}})$ we have
\begin{equation}\label{eq 3.7.2}
 (A_{21})^2+ A_{22}A_{23}+A_{23}A_{22}=0.
 \end{equation}
Now, equations (\ref{eq 3.7.1}) and (\ref{eq 3.7.2}) give us  $A_{21}^*= (-A_{21})^3$.\\
 On the other hand, using the condition  $ \alpha(\lambda_{a_{2}}) \alpha(\lambda_{a_{2}^{-1}})=\lambda_e \ot 1_{\mathbb{Q}}$, we deduce that
 \begin{equation}\label{eq 3.7.3}
  A_{23}A_{21}^*=A_{22}A_{21}^*=0.
  \end{equation}
   Hence, 
   \begin{eqnarray*}
   A_{21}^*A_{21}^* &=& A_{21}(A_{22}A_{23}+A_{23}A_{22})A_{21}^* \ (using \ (\ref{eq 3.7.1}))\\
   &=& (A_{21}A_{22})(A_{23}A_{21}^*)+ (A_{21}A_{23})(A_{22}A_{21}^*)\\
   &=& 0 \ (by \ (\ref{eq 3.7.3})). 
   \end{eqnarray*}
   

   Thus, the equality $A_{21}^*= (-A_{21})^3 $ implies $  A_{21}^*=A_{21}=0$ and also we have $A_{12}=A_{13}=0$ applying $\kappa$. This reduces the fundamental unitary to \\
 $$ 
 \begin{pmatrix}
 A_{11} & 0 & 0  & \cdots & A_{1(2l-2)} & A_{1(2l-1)} \\
 0 & A_{22} & A_{23}  & \cdots & A_{2(2l-2)} & A_{2(2l-1)}\\
 0 & A_{23}^* & A_{22}^*  & \cdots & A_{2(2l-1)}^* & A_{2(2l-2)}^*\\
 \vdots  && \hspace{1cm}        \vdots \\
 A_{l1} & A_{l2} & A_{l3} &  \cdots & A_{l(2l-2)} & A_{l(2l-1)}\\
 A_{l1}^* & A_{l3}^* & A_{l2}^*  & \cdots &  A_{l(2l-1)}^* & A_{l(2l-2)}^*\\ 
 \end{pmatrix}.$$\\\\
 The rest of the proof will be similar to Theorems \ref{thm 3 for free} and \ref{thm 2 for free}, hence omitted. \qed\\
 \brmrk
 The above results allow us to compute the cases $\Gamma=\underbrace{(\mathbb{Z}_{n_1}\ast\mathbb{Z}_{n_1}\cdot\cdot\ast\mathbb{Z}_{n_1})}_{k_1 \ copies}$ $\ast$ $\underbrace{(\mathbb{Z}_{n_2}\ast\mathbb{Z}_{n_2}\cdot\cdot\ast\mathbb{Z}_{n_2})}_{k_2 \ copies} \ast\cdot\cdot\cdot\ast \underbrace{(\mathbb{Z}_{n_l}\ast\mathbb{Z}_{n_l}\cdot\cdot\ast\mathbb{Z}_{n_l})}_{k_l \ copies}  $ where $n_1\neq n_2\neq\cdot\cdot \neq n_l$ except the case $n_1=2, \ n_l=\infty$. From \cite{free cyclic} we also came to know about $\mathbb{Q}\underbrace{(\mathbb{Z}_n \ast \mathbb{Z}_n \ast \cdot\cdot \mathbb{Z}_n)}_{k \ copies}$ except the case $n=4$. We discuss the $n=4$ case in \cite{mandal 2}.  
 \ermrk

\section{Computations for direct product of cyclic groups}\label{direct sec}
Let us make the convention of calling $\mathbb{Q}(\Gamma)$ commutative or non-commutative if the underlying $C^*$-algebra is commutative or non-commutative respectively.
In this section we will give a necessary and sufficient condition on $\Gamma$ such that $\mathbb{Q}(\Gamma)$  is  commutative.\\
\bthm \label{thm 1 for dir}
If \   $\Gamma=\underbrace{(\mathbb{Z}_n \times \mathbb{Z}_n \times \cdot\cdot\cdot \mathbb{Z}_n)}_{k \ copies}$ where $n \neq 2,4$ then $\mathbb{Q}(\Gamma) \cong  C(ISO(\Gamma))$.
\ethm
 {\it Proof:}\\
 For simplicity we present the case $k=2$. Let $\Gamma=<a,b>$ with $o(a)=o(b)=n$. The general case will follow by using similar arguments. Write the fundamental unitary as \\
  $$ U = \begin{pmatrix}
A & B & C & D\\
B^* & A^* & D^* & C^*\\
E & F & G & H\\
F^* & E^* & H^* & G^*\\
\end{pmatrix} = ((u_{ij})).$$\\
We need a few lemmas to proceed further.\\
\blmma\label{same row 0}
$AB=BA=CD=DC=EF=FE=GH=HG=0$.
\elmma
{\it Proof:}\\
Using the relation $\alpha(\lambda_{a})\alpha(\lambda_{a^{-1}})=\alpha(\lambda_{a^{-1}})\alpha(\lambda_{a})=\lambda_e \ot 1_\mathbb{Q}$ and comparing the coefficients of $\lambda_{a^{2}},\lambda_{a^{-2}},\lambda_{b^{2}},\lambda_{b^{-2}}$ we have 
$$AB^*=BA^*=B^*A=A^*B=CD^*=DC^*=C^*D=D^*C=0.$$
Applying the antipode on these relations one can deduce
  $$AB=BA=EF=FE=0.$$\\
 By similar arguments using the condition $\alpha(\lambda_{b})\alpha(\lambda_{b^{-1}})=\alpha(\lambda_{b^{-1}})\alpha(\lambda_{b})$ $=\lambda_e \ot 1_\mathbb{Q}$ we get $CD=DC=GH=HG=0$.\ \qed \\
\blmma\label{product lemma}
Product of any two different elements of each row and column of the unitary $U$ is zero. i.e. $u_{ij}u_{ik}=0 \ \forall \ j\neq k$ and $u_{ji}u_{ki}=0 \ \forall \ j\neq k$.
\elmma
{\it Proof:}\\
From the condition $\alpha(\lambda_{ab})=\alpha(\lambda_{ba})$ one can deduce 
$AE=EA$, comparing the coefficient of $\lambda_{a^{2}}$ on both sides. Hence $C^*A^*=A^*C^*$, which shows that $ AC=CA $ (applying $\kappa$ and then taking adjoint). Further, using  $\alpha(\lambda_{a})\alpha(\lambda_{a^{-1}})= \lambda_e \ot 1_\mathbb{Q}$  comparing the coefficient of $\lambda_{ab^{-1}}$ on both sides one can get $AC^*+DB^*=0$. We have
\begin{eqnarray*} 
AC^*A^* &=& -DB^*A^* \\
&=& -D(AB)^*\\
&=& 0 \ (by \ Lemma \ \ref{same row 0}).
\end{eqnarray*}
Now,
 \begin{eqnarray*}
 (AC)(AC)^* &=& (AC)(C^*A^*)\\
 &=& (CA)(C^*A^*) \ (as \ AC=CA)\\
 &=&  0 \ (as \ AC^*A^*=0).
 \end{eqnarray*}
  Hence, $AC=CA=0$ and $AE=EA=0$ (taking $\kappa$ and $\ast$ ).\\
Moreover, comparing the coefficient of $\lambda_{a^{2}}$ in $\alpha(\lambda_{ab^{-1}})=\alpha(\lambda_{b^{-1}a})$ we get $AF^*=F^*A $ and $ AD=DA $ applying the antipode and adjoint respectively. We next obtain $AD^*+CB^*=0$ comparing the coefficient of $\lambda_{ab}$ on both sides of the equation  $\alpha(\lambda_{a})\alpha(\lambda_{a^{-1}})= \lambda_e \ot 1_\mathbb{Q}$, which implies $ AD^*A^*=0 $ as $ B^*A^*=0$. Furthermore,
\begin{eqnarray*}
 (AD)(AD)^* &=& (AD)(D^*A^*)\\
 &=& (DA)(D^*A^*) \ (as \ AD=DA)\\
 &=&  0 \ (as \ AD^*A^*=0).
 \end{eqnarray*}
 Thus $AD=DA=0$ and $AF^*=F^*A=0$ (taking $\kappa$ and $\ast$ ).\\
Next, we want to show that $BD=DB=BC=CB=0$. From  $\alpha(\lambda_{ab})=\alpha(\lambda_{ba})$ and comparing the coefficient of $\lambda_{a^{-2}}$ we deduce $BF=FB$, which gives $DB=BD $ (taking $\kappa$)
and $BD^*+CA^*=0 $, implying $BD^*B^*=0$. Then $(BD)(BD)^*=D(BD^*B^*)=0$ as $BD=DB$ and $ BD^*B^*=0$. So we  get $BD=DB=0$. Similarly, we have $BC=CB=0$.\\
We also obtain similar relations replacing $A,B,C,D$ by $E,F,G,H$ respectively. \ \qed.\\
\blmma\label{normal par lem}
All the entries of $U$ are normal and partial isometries. i.e.  $u_{ij}u_{ij}^*=u_{ij}^*u_{ij}, \  u_{ij}u_{ij}^*u_{ij}=u_{ij} \ \forall \ i,j$. 
\elmma
{\it Proof:}\\
The unitarity of $U$ gives us $$AA^*+BB^*+CC^*+DD^*=1,$$  
$$A^*A+B^*B+C^*C+D^*D=1.$$
 Using Lemma \ref{product lemma} we have $ A^*A^2=A$ and $A(A^*)^2=A^*$. Now $AA^*= A^*A^2A^*=A^*A $, i.e. $A$ is a normal element. We can prove normality of other elements by similar arguments. Now we are going to show that all of them are partial isometries. First note that  $AB^*=0$ by Lemma \ref{same row 0} and Proposition \ref{normal comp}.  We claim $AC^*=AD^*=0$ too, which will imply $AA^*A=A$, multiplying by A on the left side of equation  $A^*A+B^*B+C^*C+D^*D=1$. Moreover, 
 \begin{eqnarray*}
 (AC^*)(AC^*)^* &=& AC^*CA^*\\
 &=& ACC^*A^* \ (as \ CC^*=C^*C)\\
 &=& 0 \ (as \ AC=0 \ by \ Lemma \ \ref{product lemma}).
 \end{eqnarray*}
 This shows that  $ AC^*=0$. By the same argument we can deduce $AD^*=0$.\\
Similar arguments will work for the other elements. \qed\\
\blmma
A and B commute with each elements of the set $\lbrace G,H,G^*,H^* \rbrace$. Similarly, E and F commute with each elements of the set $\lbrace C,D,C^*,D^* \rbrace$. 
\elmma
{\it Proof:}\\
From the relation  $\alpha(\lambda_{ab})=\alpha(\lambda_{ba})$ comparing the coefficient of $\lambda_{ab} $ on both sides we obtain
\begin{equation}\label{eq 4.5 1}
AG+CE=GA+EC.
\end{equation}
This gives $ A^2G+ACE=AGA+AEC$, multiplying by A on the left side of (\ref{eq 4.5 1}), hence $   A^2G=AGA $ as $AC=AE=0$ by Lemma \ref{product lemma}. \\ 
On the other hand, multiplying by A on the right side of (\ref{eq 4.5 1}) we get $AGA=GA^2$. Thus, G commutes with $A^2$ and 
\begin{equation}\label{eq 4.5 2}
A^2G=GA^2=AGA.
\end{equation}
 Moreover, G commutes with $(A^*)^2$ also by Lemma \ref{normal par lem}.  
 Now,
\begin{eqnarray*}
AG &=& A^*A^2G \ (using \ Lemma \ \ref{normal par lem}) \\
&=& A^*AGA \ (by \ (\ref{eq 4.5 2}) ) \\
&=& A^2(A^*)^2GA \ (by \ Lemma \ \ref{normal par lem}) \\
&=& GA^2(A^*)^2A \ (as \ G \ commutes \ with \ A^2 \ and \ (A^*)^2) \\
&=& GAA^*A \ (using \ Lemma \ \ref{normal par lem})\\
&=& GA.
\end{eqnarray*}
  We obtain the remaining commutation relations in a similar way. \qed \\
By the above lemmas together with Proposition \ref{iso prop}, the proof of Theorem \ref{thm 1 for dir} is completed. \qed
\brmrk
From the proof it is easily seen that we can take also $n=\infty$. For finite $n$ we get one extra relation $u_{ij}^*=u_{ij}^{n-1}$, where $((u_{ij}))$ denotes the fundamental unitary.  A close look at the proof will reveal that the fact $a^2,a^{-2},b^2,b^{-2}$ are different elements plays a crucial role here. For this reason it does not work for the cases $n=2,4$.
\ermrk
 \bthm \label{thm 2 for dir}
 Let $\Gamma=\Gamma_1 \times \Gamma_2 \cdot\cdot\cdot \times \Gamma_l  $ where  $\Gamma_i=\underbrace{(\mathbb{Z}_{n_i}\times\mathbb{Z}_{n_i}\cdot\cdot\times\mathbb{Z}_{n_1})}_{k_i \ copies}$. Also assume that  $n_1\neq n_2\neq\cdot\cdot \neq n_l$ and $n_i\neq2,4 \ \forall \ i$, then $\mathbb{Q}(\Gamma)$ will be $C(ISO\underbrace{(\mathbb{Z}_{n_1}\times\mathbb{Z}_{n_1}\cdot\cdot\times\mathbb{Z}_{n_1})}_{k_1 \ copies})$ $\hat{\ot} C(ISO \underbrace{(\mathbb{Z}_{n_2}\times\mathbb{Z}_{n_2}\cdot\cdot\times\mathbb{Z}_{n_2})}_{k_2 \ copies}) \hat{\ot} \cdot\cdot \hat{\ot} C(ISO \underbrace{(\mathbb{Z}_{n_l}\times\mathbb{Z}_{n_l}\cdot\cdot\times\mathbb{Z}_{n_l})}_{k_l \ copies} )  $.  i.e. $\mathbb{Q}(\Gamma) \cong \mathbb{Q}(\Gamma_1) \hat{\ot} \mathbb{Q}(\Gamma_2)\hat{\ot} \cdot\cdot\cdot \hat{\ot} \mathbb{Q}(\Gamma_l)$.    
 \ethm
 {\it Proof:}\\
As before we give the proof for the case where all $k_i=1$ to simplify notation. Now $\Gamma=\mathbb{Z}_{n_1}\times \mathbb{Z}_{n_{2}}\times \cdot\cdot\cdot \times \mathbb{Z}_{n_l}$ and choose $\lbrace a_1,a_1^{-1},a_2,a_2^{-1},\cdot\cdot\cdot,a_l,a_l^{-1}\rbrace$ as the standard symmetric generating set of $\Gamma$, where $o(a_i)=n_i \ \forall \ i$. \\
Fundamental unitary is of the form \\
  $$ U=
 \begin{pmatrix}
 A_{11} & A_{12} & A_{13} & A_{14} & \cdots & A_{1(2l-1)} & A_{1(2l)}\\
 A_{12}^* & A_{11}^* & A_{14}^* & A_{13}^* & \cdots &  A_{1(2l)}^* & A_{1(2l-1)}^*\\
 A_{21} & A_{22} & A_{23} & A_{24} & \cdots & A_{2(2l-1)} & A_{2(2l)}\\
 A_{22}^* & A_{21}^* & A_{24}^* & A_{23}^* & \cdots & A_{2(2l)}^* & A_{2(2l-1)}^*\\
 \vdots  && \hspace{1cm}        \vdots \\
 A_{l1} & A_{l2} & A_{l3} & A_{l4} & \cdots & A_{l(2l-1)} & A_{l(2l)}\\
 A_{l2}^* & A_{l1}^* & A_{l4}^* & A_{l3}^* & \cdots &  A_{l(2l)}^* & A_{l(2l-1)}^*\\ 
 \end{pmatrix}.$$\\\\
First we want to reduce it to a block diagonal form. By Remark \ref{rmk in sec 1} this will complete the proof, because in this case $\mathbb{Q}(\Gamma_i)$ is commutative by Theorem \ref{thm 1 for dir}, hence $\ot ^{max}$ coincides with $\hat{\ot}$.\\
 We remark that
\begin{equation}\label{eq 4.7.1}
A_{1i}A_{1j}=A_{1j}A_{1i}=0 \ \forall \ i\neq j
\end{equation}
by arguments similar to those in the proof of Lemma \ref{product lemma}. Using (\ref{eq 4.7.1}) we have   
 $\alpha(\lambda_{a_{1}^{n_1-1}})= \lambda_{a_{1}^{n_1-1}}\ot A_{11}^{n_1-1}+ \lambda_{({a_{1}^{-1}})^{n_1-1}}\ot A_{12}^{n_1-1} +\cdot\cdot\cdot+ \lambda_{a_{l}^{n_1-1}}\ot A_{1(2l-1)}^{n_1-1}+ \lambda_{({a_{l}^{-1}})^{n_1-1}}\ot A_{1(2l)}^{n_1-1}$. Now, from the relation $\alpha(\lambda_{a_{1}^{n_1-1}})=\alpha(\lambda_{a_{1}^{-1}})$ one obtains
 \begin{equation}\label{eq 4.7.2}
  A_{12}^*=A_{12}^{n_1-1}, \ A_{11}^*=A_{11}^{n_1-1}, \ A_{1i}^*=0 \ \forall \ i>2
  \end{equation}
  by comparing the coefficients of $\lambda_{a_{1}},\lambda_{a_{1}^{-1}},\lambda_{a_{2}},\lambda_{a_{2}^{-1}},\cdot\cdot\cdot,\lambda_{a_{l}},\lambda_{a_{l}^{-1}}$. Applying the antipode on (\ref{eq 4.7.2}) we find $A_{i1}=A_{i2}=0 \ \forall \ i >1$. Now the fundamental unitary reduces to \\
$$ \begin{pmatrix}
 A_{11} & A_{12} & 0 & 0 & \cdots & 0 & 0\\
 A_{12}^* & A_{11}^* & 0 & 0 & \cdots &  0 & 0\\
 0 & 0 & A_{23} & A_{24} & \cdots & A_{2(2l-1)} & A_{2(2l)}\\
 0 & 0 & A_{24}^* & A_{23}^* & \cdots & A_{2(2l)}^* & A_{2(2l-1)}^*\\
 \vdots  && \hspace{1cm}        \vdots \\
 0 & 0 & A_{l3} & A_{l4} & \cdots & A_{l(2l-1)} & A_{l(2l)}\\
 0 & 0 & A_{l4}^* & A_{l3}^* & \cdots &  A_{l(2l)}^* & A_{l(2l-1)}^*\\ 
 \end{pmatrix}.$$\\\\
Proceeding in a similar way we get the desired block diagonal form. Note that the last relation $\alpha(\lambda_{a_{l}^{n_l-1}})=\alpha(\lambda_{a_{l}^{-1}})$ has not been used. For this reason we can also include the case $n_l=\infty$.   \qed \\ 
\bthm \label{thm 3 for dir}
The conclusion of Theorem \ref{thm 2 for dir} remains valid if $n_1=2$ where $2<n_2<n_3\cdot\cdot<n_l \leq \infty$ and  $n_i\neq 4 \ \forall \ i$.
\ethm
 {\it Proof:}\\
 For simplicity consider $k_i=1 \ \forall \ i$ so that $\Gamma=\mathbb{Z}_2\times \mathbb{Z}_{n_2}\times \mathbb{Z}_{n_3}\times \cdot\cdot\cdot \times \mathbb{Z}_{n_l}$ and $o(a_1)=2$.\\
We only show how to reduce the fundamental unitary to the form given below, rest of the arguments are similar to those of Theorem \ref{thm 2 for dir}. \\
 $$ \begin{pmatrix}
 A_{11} & 0 & 0  & \cdots & 0 & 0 \\
 0 & A_{22} & A_{23}  & \cdots & 0 & 0\\
 0 & A_{23}^* & A_{22}^*  & \cdots & 0 & 0\\
 \vdots  && \hspace{1cm}        \vdots \\
 0 & 0 & 0 &  \cdots & A_{l(2l-2)} & A_{l(2l-1)}\\
 0 & 0 & 0  & \cdots &  A_{l(2l-1)}^* & A_{l(2l-2)}^*\\ 
 \end{pmatrix}$$\\\\
 Note that product of any two different elements of each row is zero (except possibly the first row) following the lines of arguments of Lemma \ref{product lemma}, which means
 \begin{equation}\label{eq 4.8.1}
 A_{k1}A_{ki}=0 \ \forall \ k, i > 1.
 \end{equation}
 Furthermore, considering the coefficient of $\lambda_e$ in the expression of $\alpha(\lambda_{a_{k}^2})$ we obtain
 \begin{equation}\label{eq 4.8.2}
 A_{k1}^2=0 \ \forall \ k >1.  
 \end{equation}
 Our aim is to show that $A_{k1}=0 \ \forall \ k>1$. Now, using the condition $ \alpha(\lambda_{a_{k}}) \alpha(\lambda_{a_{k}^{-1}})=\lambda_e \ot 1_{\mathbb{Q}} \ \forall \ k>1$, one can deduce
 \begin{equation}\label{eq 4.8.3}
  A_{k1}A_{k1}^*+A_{k2}A_{k2}^*+\cdot\cdot\cdot + A_{k(2l-1)}A_{k(2l-1)}^*=1.
  \end{equation}
  Multiplying by $A_{k1}$ on the left side of (\ref{eq 4.8.3}) we get 
  \begin{equation}\label{eq 4.8.4}
   A_{k1}(A_{k1}A_{k1}^*+A_{k2}A_{k2}^*+\cdot\cdot\cdot + A_{k(2k-1)}A_{k(2k-1)}^*)=A_{k1} \ \forall \ k>1.
  \end{equation}
  Now using (\ref{eq 4.8.1}) and (\ref{eq 4.8.2}) one can find $A_{k1}=A_{k1}^*=0 \ \forall \ k>1 $, hence $A_{12}=A_{13}=A_{14}=A_{15}=\cdot\cdot A_{1(2l-2)}=A_{1(2l-1)}=0$ by applying the antipode.\\
This gives one-step reduction of the fundamental unitary to the following form \\
 $$ U=
 \begin{pmatrix}
 A_{11} & 0 & 0  & \cdots & 0 & 0 \\
 0 & A_{22} & A_{23}  & \cdots & A_{2(2l-2)} & A_{2(2l-1)}\\
 0 & A_{23}^* & A_{22}^*  & \cdots & A_{2(2l-1)}^* & A_{2(2l-2)}^*\\
 \vdots  && \hspace{1cm}        \vdots \\
 0 & A_{l2} & A_{l3} &  \cdots & A_{l(2l-2)} & A_{l(2l-1)}\\
 0 & A_{l3}^* & A_{l2}^*  & \cdots &  A_{l(2l-1)}^* & A_{l(2l-2)}^*\\ 
 \end{pmatrix}.$$\\\\
 Then we proceed similarly as Theorem \ref{thm 2 for dir} to achieve the desired reduction .   \qed\\
 \brmrk
 Notice that unlike the free case, the above proof even works for the case when $o(a_l)=\infty$.
 \ermrk
\bthm
The conclusion of Theorem \ref{thm 2 for dir} is valid for $n_1=4$  where $4\neq n_2\neq\cdot\cdot \neq n_l$ and $n_i\neq2 \ \forall \ i$, if we replace $\hat{\ot}$ by $\ot^{max}$, i.e. $\mathbb{Q}(\Gamma) \cong \mathbb{Q}(\Gamma_1) \ot^{max} \mathbb{Q}(\Gamma_2) \ot^{max}\cdot\cdot\cdot \ot^{max} \mathbb{Q}(\Gamma_l)$. However, in this case $\mathbb{Q}(\Gamma)$ is non-commutative.
\ethm
 {\it Proof:}\\
 As before, assume without loss of generality and for the sake of simplicity of exposition that $k_i=1 \ \forall \ i $ and $o(a_1)=4$. 
 Again, the product of any two different elements of each row and column is zero (except the first two rows and columns) by the arguments of Lemma \ref{product lemma}. Hence all entries except possibly the first two rows and columns are normal by arguments similar to those of Lemma \ref{normal par lem}. Now we claim that $A_{k1}=A_{k2}=0 \ \forall \ k>1$. Considering the coefficients of $\lambda_{a_1}, \lambda_{a_{1}^{-1}}$ of $\alpha(\lambda_{a_{k}^3}) \ \forall \ k>1$ we have
 \begin{equation}\label{eq 4.10.1}
  (A_{k1}^2+A_{k2}^2)A_{k1} = (A_{k1}^2+A_{k2}^2)A_{k2}=0.
  \end{equation}
  This gives us $ A_{k1}^3=A_{k2}^3=0 $, hence $ A_{k1}=A_{k2}=0 \ \forall \ k>1$, as they are normal. Applying the antipode we deduce $A_{1k}=0 \ \forall \ k>2$.\\
 Rest of the proof is very similar to Theorem \ref{thm 2 for dir}, hence omitted.  \qed\\
 
Combining the above theorems in this section and Propositions \ref{prop 1} and \ref{prop 4} we get the following necessary and sufficient condition for $\mathbb{Q}(\Gamma)$ to be commutative.
 \bcrlre 
   $\mathbb{Q}(\Gamma)$ is commutative if and only if $\Gamma$ must be of the form $\underbrace{(\mathbb{Z}_{n_1}\times\mathbb{Z}_{n_1}\cdot\cdot\times\mathbb{Z}_{n_1})}_{k_1 \ copies}$ $\times$ $\underbrace{(\mathbb{Z}_{n_2}\times\mathbb{Z}_{n_2}\cdot\cdot\times\mathbb{Z}_{n_2})}_{k_2 \ copies} \times\cdot\cdot\cdot\times \underbrace{(\mathbb{Z}_{n_l}\times\mathbb{Z}_{n_l}\cdot\cdot\times\mathbb{Z}_{n_l})}_{k_l \ copies}  $ where $n_i \neq 4 \ \forall \ i$ and if $n_j=2$ for some $j$, then $k_j$ must be $1$.
 \ecrlre

\section {Examples of $(\Gamma,S)$ for which $\mathbb{Q}(\Gamma)\cong\mathcal{D}_{\theta}(C^*(\Gamma), \Delta_{\Gamma})$}
It has been observed in Proposition 2.3 of \cite{doubling} that, if there exists a non trivial automorphism of order $2$ which preserves the generating set, then $\mathcal{D}_{\theta}(C^*(\Gamma), \Delta_{\Gamma})$ is always a quantum subgroup of $\mathbb{Q}(\Gamma)$. For many examples studied by the authors of \cite{grp algebra}, \cite{S_n}, \cite{dihedral} $\mathbb{Q}(\Gamma)$ coincides with the doubled group algebra. In this section we produce more examples of groups where this occurs.\\
\subsection{Dihedral groups with two different generating sets}
Dihedral group has two presentations\\
\begin{eqnarray}
D_{2n}= <a,b| \ a^2=b^n=e,ab=b^{-1}a>, \label{di 1} \\ 
D_{2n}= <s,t| \ s^2=t^2=(st)^n=e>, \label{di 2}
\end{eqnarray}
where $e$ denotes the identity element of the group. In \cite{dihedral} the authors calculated QISO for $D_{2(2n+1)}$ with the presentation (\ref{di 1}). Let us calculate it for $D_{2n}$ with presentation (\ref{di 2}).\\
\bthm
Let $D_{2n}=<s,t| \ s^2=t^2=(st)^n=e >$, then its QISO is isomorphic to $\mathcal{D}_{\theta}(C^*(D_{2n}), \Delta_{D_{2n}})$ with respect to the automorphism $\theta$ given by $\theta(s)=t, \  \theta(t)=s$. 
\ethm
{\it Proof:}\\
The action is defined by\\
$$\alpha(\lambda_{s})= \lambda_{s} \ot A + \lambda_{t} \ot B,$$
$$\alpha(\lambda_{t})= \lambda_{s} \ot C + \lambda_{t} \ot D.$$ \\
 Here $\begin{pmatrix}
A & B  \\
C & D \\
\end{pmatrix}$ 
is the corresponding fundamental unitary.\\
We present the case for odd $n=2k+1$, as the proof is almost the same for even n.\\
We begin with the simple observation that given $s^2=t^2=e$, the condition $(st)^n=e$ is equivalent to $(st)^{k}s=(ts)^{k}t$.
From the relation $\alpha(\lambda_{s})=\alpha(\lambda_{s^{-1}})$ we have $A=A^*, \ B=B^*$ and using $\alpha(\lambda_{s^{2}})=\lambda_e \ot 1_{\mathbb{Q}}$ we get $A^2+B^2=1, \ AB=BA=0$.\\
Similarly $C=C^*, \ D=D^*$ and $C^2+D^2=1, \ CD=DC=0$.\\
Applying the antipode on the above equations we find that $$A^2+C^2=B^2+D^2=1,$$ $$ AC=CA=BD=DB=0.$$
So we obtain $A^2=D^2, \ B^2=C^2$ and clearly $A^2,B^2,C^2,D^2$ are central projections.
Now we are going to use $(st)^n=(st)^{2k+1}=e$, from which we deduce
\begin{equation} \label{di 3}
(st)^ks=(ts)^kt.
\end{equation} 
We want to obtain analogues of (\ref{di 3}) with $(s,t)$ replaced by $(A,D)$ as well as $(B,C)$. Using relation (\ref{di 3}) we have
\begin{equation} \label{di 4} 
(AD)^{k}A+(BC)^{k}B=(DA)^{k}D + (CB)^{k}C.
\end{equation}
Applying $\kappa$ on (\ref{di 4}) we get
\begin{equation} \label{di 5} 
  (AD)^{k}A+(CB)^{k}C=(DA)^{k}D+(BC)^{k}B.
 \end{equation}
Equations (\ref{di 4}) and (\ref{di 5}) together imply $(AD)^{k}A=(DA)^{k}D, \ (CB)^{k}C=(BC)^{k}B$.\\
Now it follows from Lemma \ref{doub lemma} that $\mathbb{Q}(D_{2n})$ coincides with  $\mathcal{D}_{\theta}(C^*(D_{2n}),\Delta_{D_{2n}})$  corresponding to the order 2 automorphism $\theta$ given by $\theta(s)=t, \ \theta(t)=s$. \qed\\
\brmrk
We get the same result for $D_{2n}$ with presentation (\ref{di 1}) except $n=4$ case. The proof is very similar to the case $\mathbb{Z}_2 \times \mathbb{Z}_n$, hence omitted. This extends the result of \cite{dihedral}.
\ermrk
\subsection{Baumslag-Solitar group}
The group has the presentation $\Gamma=<a,b| \ b^{-1}ab=a^{2}>$.\\
We can easily get the following relations among the generators, \\
$$b^{-1}ab=a^{2}, \  b^{-1}a^{-1}b=a^{-2},$$
$$ab=ba^{2}, \ b^{-1}a=a^{2}b^{-1}, \ a^{-1}b^{-1}a=ab^{-1},$$
$$b^{-1}a^{-1}=a^{-2}b^{-1}, \ a^{-1}b=ba^{-2}, \ ab=ba^{2}, \ aba^{-1}=ba.$$
Write the fundamental unitary as 
$$\begin{pmatrix}
A & B & C & D\\
B^* & A^* & D^* & C^*\\
E & F & G & H\\
F^* & E^* & H^* & G^*\\
\end{pmatrix}.$$\\
Our aim is to show $C=D=E=F=0$.\\
Using the relation $\alpha(\lambda_{ab})=\alpha(\lambda_{ba^2})$ and comparing the coefficients of   $\lambda_{b^{2}}$ and $\lambda_{b^{-2}}$, we obtain $CG=DH=0 $ and by applying the antipode we also get $ EG=HF=0 $.\\
Further, using the condition $\alpha(\lambda_{a^{-1}b^{-1}a})=\alpha(\lambda_{ab^{-1}})$ one can find $CH^*=DG^*=0,$ hence $ EH=GF=0$. Also note that $EF=0$ arguing along the lines of Lemma \ref{lemma for free}. So $E(GG^*+HH^*+FF^*+EE^*)=E $, implying $ E^2E^*=E$  as $(GG^*+HH^*+FF^*+EE^*)=1, \ EG=EF=EH=0$. Similarly, we obtain $F^*F^2=F$.  \\
Next, we compare the coefficients of $\lambda_{b^{2}}$ and  $\lambda_{b^{-2}}$ on both sides of  $\alpha(\lambda_{b^{-1}ab})=\alpha(\lambda_{a^{2}})$ to have $C^2=D^2=0 $, which implies $ E^2=F^2=0$ using antipode and taking adjoint of the elements.\\
Now using the above relations, $C=D=E=F=0$ as $ E^2E^*=E, \ F^*F^2=F$.\\
 The fundamental unitary reduces to the form 
$$\begin{pmatrix}
A & B & 0 & 0\\
B^* & A^* & 0 & 0\\
0 & 0 & G & H\\
0 & 0 & H^* & G^*\\
\end{pmatrix}.$$\\
Moreover, we claim that $H=0$. Using  $\alpha(\lambda_{ab})=\alpha(\lambda_{ba^2})$ we have\\
$AH=BH=0$, equating the coefficients of $\lambda_{ab},\lambda_{a^{-1}b^{-1}}$ on both sides. Thus
$(A^*A+B^*B)H=H$, that is, $ H=0$.\\
This gives the following reduction: 
$$\begin{pmatrix}
A & B & 0 & 0\\
B^* & A^* & 0 & 0\\
0 & 0 & G & 0\\
0 & 0 & 0 & G^*\\
\end{pmatrix}.$$\\
 Moreover, using the relations among the generators we deduce that $G^*AG=A^2, \ G^*BG=B^2, \ AG=GA^2, \ BG=GB^2$ and also that $AA^*, \ BB^*$ are central projections of the algebra.\\
It now follows from Lemma \ref{doub lemma} that $\mathbb{Q}(\Gamma)$ is isomorphic to $\mathcal{D}_{\theta}(C^*(\Gamma),\Delta_{\Gamma})$ with respect to the automorphism $\theta$ given by $a\mapsto a^{-1}, \  b\mapsto b$. \qed\\ 
\subsection{Some groups of the form $<a,b| \ o(a)=2,o(b)=3>$} \label{tetra}
First we conclude a lemma which will be useful for the proof of Theorem \ref{tet,oct th}.\\
\blmma\label{order 3}
If \  $\Gamma=<a,b| \ o(a)=2,o(b)=3>$, then its fundamental unitary must be of the form\\
$$\begin{pmatrix}
A & 0 & 0 \\
0 & E & F \\
0 & F^* & E^* \\
\end{pmatrix}.$$
\elmma
{\it Proof:}\\
Using the relation $\alpha(\lambda_{b^{2}})=\alpha(\lambda_{b^{-1}})$, and
comparing the coefficients of $\lambda_a$ and $\lambda_{a^{-1}}$ from both sides we will get the reduced block diagonal form. \qed\\
Now we use the above lemma to compute quantum isometry groups of some concrete examples.
\bthm\label{tet,oct th}
Let $\Gamma$ be as in the statement of Lemma \ref{order 3} and consider the automorphism $\theta$ defined by $\theta(a)=a, \ \theta(b)=b^{-1}$. Then $\mathbb{Q}(\Gamma)\cong\mathcal{D}_{\theta}(C^*(\Gamma), \Delta_{\Gamma})$, for the following three examples:\\
$1. (ab)^3=1$, (Tetrahedral)\\
$2. (ab)^4=1$, (Octahedral)\\
$3. (ab)^5=1$. (Icosahedral)\\
\ethm
{\it Proof:}\\
In each of these cases, we can apply Lemma \ref{order 3} to get $A,E,F$ such that the action $\alpha$ is given by\\
$$\alpha(\lambda_{a})= \lambda_{a} \ot A , $$
$$\alpha(\lambda_{b})=  \lambda_{b} \ot E + \lambda_{{b}^{-1}} \ot F ,$$
$$\alpha(\lambda_{{b}^{-1}})=  \lambda_{b} \ot F^* + \lambda_{{b}^{-1}} \ot E^* .$$\\

 
Also $A^2=1, \ A=A^*$, applying $\alpha(\lambda_{a})=\alpha(\lambda_{a^{-1}})$ and $\alpha(\lambda_{a^{2}})=\lambda_e \ot 1_{\mathbb{Q}}$.
Similarly $E^2=E^*, \ F^2=F^*$ using the condition $\alpha(\lambda_{b^{2}})=\alpha(\lambda_{b^{-1}})$. 
We also get $EF=FE=0$ arguing along the lines of Lemma \ref{same row 0}. 

Now consider the relation $(ab)^n=1 $ ($n=3,4,5$ respectively), which gives $ ab=(b^{-1}a)^{n-1}$ and  $(ba)=(ab^{-1})^{n-1}$.
Using these relations we can deduce
\begin{equation}\label{eq 5.3.1}
AE= \underbrace{(E^*A)(E^*A)\cdot\cdot\cdot(E^*A)}_{(n-1) \ times},
\end{equation}
\begin{equation}\label{eq 5.3.2}
 EA=\underbrace{(AE^*)(AE^*)\cdot\cdot\cdot(AE^*)}_{(n-1) \ times}.
 \end{equation}
 Moreover,
 \begin{eqnarray*}
  A(EE^*) &=& \underbrace{(E^*A)(E^*A)\cdot\cdot\cdot(E^*A)}_{(n-1) \ times}E^* \ (by \ (\ref{eq 5.3.1}))\\
 &=& E^*\underbrace{(AE^*)(AE^*)\cdot\cdot\cdot(AE^*)}_{(n-1) \ times}\\
 &=& E^*(EA) \ (using \ (\ref{eq 5.3.2}))\\
 &=& (EE^*)A \ (as \ E^*=E^2). 
\end{eqnarray*}  
Hence, $EE^*$ is a central projection. By similar arguments it can be proved that $FF^*$ is a central projection. By Lemma \ref{doub lemma} we get the isomorphism between $\mathbb{Q}(\Gamma)$ and $\mathcal{D}_{\theta}(C^*(\Gamma),\Delta_{\Gamma})$  with respect to the automorphism $\theta$ given by $b\mapsto b^{-1}, \ a\mapsto a$. \qed\\
\section{Coxeter groups as examples of  $\Gamma$ such that $\mathbb{Q}(\Gamma)\cong (C^*(\Gamma), \Delta_{\Gamma})$}
In this section we will compute QISO for certain classes of Coxeter groups. 
The Coxeter group with parameters $(l,m,n)$ and $l\leq m\leq n$ has the following presentation
 $$\Gamma=<a,b,c| \ o(a)=o(b)=o(c)=2,(ac)^l=(ab)^m=(bc)^n=e>.$$
Here we take one special class, namely $l=2,m=3$ and $n$ is any positive integer co-prime to 6. \\
\bthm
\label{cox}
Let $\Gamma$ be the Coxeter group with parameters $(2,3,n)$ as above. Then $\mathbb{Q}(\Gamma)\cong (C^*(\Gamma), \Delta_{\Gamma})$.
\ethm
{\it Proof:}\\
The fundamental unitary is of the form 
$$\begin{pmatrix}
A & B & C \\
D & E & F \\
G & H & K \\
\end{pmatrix}.$$\\
We divide the proof into a number of lemmas.\\
\blmma
$B=D=H=F=0$.
\elmma
{\it Proof:}\\
First note that $ \alpha(\lambda_{ac})=\alpha(\lambda_{ca})$ as $(ac)^2=a^2=c^2=e$.\\
$$\alpha(\lambda_{ac})=\lambda_{e}\ot (AG+BH+CK)+\lambda_{ac}\ot (AK+CG)+ \lambda_{ab}\ot AH + \lambda_{bc}\ot BK + \lambda_{ba}\ot BG + \lambda_{cb}\ot CH,$$
$$\alpha(\lambda_{ca})=\lambda_{e}\ot (GA+HB+KC)+\lambda_{ac}\ot (KA+GC)+ \lambda_{ab}\ot GB + \lambda_{bc}\ot HC + \lambda_{ba}\ot HA + \lambda_{cb}\ot KB.$$\\
From the above equations we have
$$AH=GB, \ BK=HC, \ HA=BG, \ CH=KB.$$ 
Now as the action is length preserving, we also get
\begin{equation}\label{eq 6.2.1}
AG+BH+CK=0.
\end{equation}
Using the condition $\alpha(\lambda_{b^{2}})=\alpha(\lambda_{e})$ one obtains
$DF+FD=0$ by comparing the coefficient of $\lambda_{ac}$ on both sides, and hence,
$HB+BH=0$ (applying $\kappa$).\\
This gives us
\begin{equation}\label{eq 6.2.2}
HBH+BH^2=0.
\end{equation} 
Again using  $\alpha(\lambda_{c^{2}})=\alpha(\lambda_{e})$ we have
$G^2+H^2+K^2=1$, hence
\begin{equation} \label{eq 6.2.3}
 BG^2+BH^2+BK^2=B.
 \end{equation}
Now (\ref{eq 6.2.1}) implies $ HAG+HBH+HCK=0$, which gives
 $ (BG)G+HBH+(BK)K=0$, as $HA=BG, \ BK=HC$. Thus we get $BG^2-BH^2+BK^2=0$ using (\ref{eq 6.2.2}). Comparing it with (\ref{eq 6.2.3}) we deduce $B=2BH^2$.\\
We will show now $ BH^2=0$. Multiplying H on the right side of (\ref{eq 6.2.1}) we find \\
$AGH+BH^2+CKH=0$, where $GH=KH=0$, using  $\alpha(\lambda_{c^{2}})=\alpha(\lambda_{e})$ and comparing the coefficients of $\lambda_{ab}, \ \lambda_{bc}$ respectively.\\
So $BH^2=0$, hence $B=0$,
and $ D=0 $ too applying the antipode.\\
This also gives $HA=BG=0, \ HC=BK=0$. Moreover, $ (A^2+C^2)=1$ by using $\alpha(\lambda_{a^{2}})=\alpha(\lambda_{e})$. We have 
\begin{eqnarray*}
H &=& H(A^2+C^2) \ (as \ (A^2+C^2)=1) \\
&=& (HA)A + (HC)C \\
&=& 0 \ (as \ HA=HC=0).
\end{eqnarray*}  
 Similarly, one can get $F=0$. \qed\\
\\
\blmma
$G=C=0$.
\elmma
{\it Proof:}\\
First we need to derive some more relations among the generators from the defining ones.

From $(ab)^3=e$ we have $ aba=bab$ as $a^2=b^2=e$.
Our aim is to show that $ECE=0$. We claim that the term $bcb$ is not equal to any of the terms $aba, abc, cba, cbc$. Clearly, $bcb\neq aba$ as we have $aba=bab$. Now $cbc\neq bcb$ as $(bc)^3\neq e$. If $(bc)^3= e$ then using the hypothesis we can obtain $bc=e$ as n is co-prime to 3 too, which implies $b=c$. Furthermore, if $bcb=abc$ which implies $(bc)^2=ab$, then we get $(bc)^6=e$. Hence, one can deduce $bc=e$. Similarly, we can argue that $bcb$ can't be equal to $cba$. Now using $\alpha(\lambda_{bab})=\alpha(\lambda_{aba})$ we obtain
\begin{eqnarray}
\alpha(\lambda_{bab}) &=& \lambda_{bab}\ot EAE + \lambda_{bcb} \ot ECE \label{eq 6.1}, \\ 
\alpha(\lambda_{aba}) &=& \lambda_{aba} \ot AEA + \lambda_{abc} \ot AEC +  \lambda_{cba} \ot CEA + \lambda_{cbc} \ot CEC. \label{eq 6.2}
\end{eqnarray}
Comparing the both sides of (\ref{eq 6.1}) and (\ref{eq 6.2})  we find
$ECE=0$, also $E^2=1$ from the condition  $\alpha(\lambda_{b^{2}})=\alpha(\lambda_{e})$. Moreover,
\begin{eqnarray*}
C &=& E^2CE^2 \ (as \ E^2=1)\\
&=& E(ECE)E \\
&=& 0 \ (as \ ECE=0). 
\end{eqnarray*}
 Applying the antipode we find $G=0$.  \qed \\\\
{\it Proof of Theorem \ref{cox}}:\\
By the above lemmas, we have reduced the fundamental unitary to the following form: 
$$\begin{pmatrix}
A & 0 & 0 \\
0 & E & 0 \\
0 & 0 & K \\
\end{pmatrix}.$$\\
So by Corollary \ref{group algebra cor}, there is a CQG  isomorphism from $(C^*(\Gamma), \Delta_{\Gamma})$ to $\mathbb{Q}(\Gamma)$ sending $\lambda_a,\lambda_b,\lambda_c$ to A, E, K respectively. \ \qed\\
\section{An excursion to QISO of compact matrix quantum groups}
In this very brief last section, we want to extend the formulation of quantum isometry group to the realm of quantum groups. Let us consider a compact matrix quantum group $(\clq, \Delta)$
 which has a finite (say $n$) dimensional unitary  fundamental representation $\pi$ with $((\pi_{ij} ))$ being the corresponding unitary in $M_n(\clq)$. Indeed, by definition,
  every irreducible representation of $(\clq, \Delta)$ is a sub-representation of tensor copies of $\pi$ and $\bar{\pi}$, so as in Subsection 2.6, we may consider a central length function $l$ which takes an irreducible say $\alpha$ to the smallest non-negative integer $k$ such that $\alpha \subset \alpha_1 \ot \alpha_2 \ot \cdot\cdot\cdot  \ot \ \alpha_k$ where each $\alpha_i$ is either $\pi$ or $\bar{\pi}$. As shown in \cite{jyotish_new}, this gives rise to a spectral triple, generalizing the construction of 
     $D_\Gamma$ for a finitely generated discrete group in Subsection 2.1. Moreover, this spectral triple satisfies the condition of Theorem \ref{existence thm}, 
     hence the quantum isometry group exists. Let us denote the quantum isometry group by $\mathbb{Q}(\hat{\clq},\{ u_{ij} \})$, where $u_{ij}$ consist of both $\pi_{ij}$'s as well as 
      $\pi_{ij}^*$'s, or simply $\mathbb{Q}(\hat{\clq})$ 
     if  the matrix elements $u_{ij}$ are understood. By Gram-Schmidt, we convert $\lbrace u_{ij} \rbrace$ to an orthogonal set, say $\lbrace u_{ij}^{\prime} \rbrace$ with respect to the Haar state of $(\clq, \Delta)$. 
       Indeed, as in the group case, the action of $\mathbb{Q}(\hat{\clq})$, say $\beta$, is determined by $q^{ij}_{kl}$ such that $$ \beta(u_{ij}^{\prime})=\sum_{kl} u_{kl}^{\prime} \ot q^{ij}_{kl}.$$
       In other words, the quantum isometry group is generated by $q^{ij}_{kl}$ subject to the relations that make it  a unitary and also make 
       the above $\beta$ a $\ast$-homomorphism from $\clq$ to $\clq \hat{\ot} \ 
       \mathbb{Q}(\hat{\clq})$.
       \brmrk\label{7th sec remrk}
        Note that $\mathbb{Q}(\hat{\clq},\{ u_{ij} \})$ is the universal $C^*$-algebra generated by $q^{ij}_{kl}$ such that $Q=((q^{ij}_{kl}))$ is unitary as well as $Q^tE\bar{Q}E^{-1}=E\bar{Q}E^{-1}Q^t=I$, where $E=((E^{ij}_{kl}))=((\langle(u_{ij}^{\prime})^*,(u_{kl}^{\prime})^*\rangle))$ and $\beta$ given above is a $C^*$-homomorphism on $\clq$ by the similar argument of Proposition \ref{proposition category}. 
       \ermrk

     
     It is interesting to note that the formulation for quantum isometry group for a compact matrix quantum group allows us to consider even the group algebras with a set of generators which are 
      not necessarily of the form $\delta_g$ for elements $g$ of the group, i.e. not necessarily group-like elements of the group $C^*$-algebra. This flexibility of choice can 
        have quite interesting implications for the resulting quantum isometry groups, as illustrated by the  example below. We consider the group $\Gamma = \mathbb{Z} \times \mathbb{Z}_2$. It has a natural set of generators consisting of group-like elements as in Theorem \ref{thm 3 for dir}, where the resulting QISO turned out to be $\mathbb{Q}(\mathbb{Z}) \hat{\ot} \mathbb{Q}(\mathbb{Z}_2)$. It can be identified as the doubling of the group algebra too.   
           However, we can also view it as a matrix quantum group with a fundamental unitary whose entries are 
           not group-like elements. More precisely, note that $C_r^*(\mathbb{Z} \times \mathbb{Z}_2)$ is isomorphic with $C(\mathbb{T}) \oplus C(\mathbb{T})$ as a $C^*$-algebra, 
           and it can be described as the  $C^*$-algebra  $C^*\lbrace \gamma, \gamma^{\prime} \rbrace$ where $\gamma$, $\gamma^\prime$ denotes the canonical generators of 
            the two copies of $C(\mathbb{T})$. The $C^*$-algebra is the universal one with two generators satisfying the following relations:
            \begin{eqnarray} \label{eq 7.1}
            \gamma.\gamma^* = \gamma^\ast.\gamma , \  \gamma^{\prime}\gamma^{\prime\ast}=\gamma^{\prime\ast}\gamma^{\prime},\\
          \gamma.\gamma^{\prime}= \gamma^{\prime}.\gamma =0, \label{eq 71} \\
          \gamma.\gamma^\ast + \gamma^{\prime}.\gamma^{\prime\ast} = 1. \label{eq 72}
          \end{eqnarray}
         From the group structure of $\Gamma$, it is easy to see that $\clf:=\{ \gamma, \gamma^*, \gamma^\prime, (\gamma^\prime)^* \}$ gives the 
          matrix coefficients of a $2$-dimensional fundamental unitary representation, not consisting of group-like elements. We have the following description of the quantum isometry group for 
           the generating set $\clf$, which is again a doubling, but not of the group algebra itself.

     \bthm\label{sec 7 th}
         $\mathbb{Q}(\widehat{C_r^*(\mathbb{Z} \times \mathbb{Z}_2)}, \clf)$  is isomorphic with a  doubling of the quantum group $\mathbb{Q}(\mathbb{Z}) \star \mathbb{Q}(\mathbb{Z})$ with respect to the order $2$ automorphism $\theta$ defined by 
   $$\gamma_1 \mapsto \gamma_1^{\prime}, \gamma_1^{\prime} \mapsto \gamma_1,\gamma_2 \
\mapsto \gamma_2^{\prime},\gamma_2^{\prime} \mapsto \gamma_2,$$
where $\gamma_1, \gamma_1^{\prime}, \gamma_2, \gamma_2^{\prime}$ generate the underlying $C^*$- algebra of the CQG $\mathbb{Q}(\mathbb{Z}) \star \mathbb{Q}(\mathbb{Z})$.    
         \ethm
         {\it Proof:}
           Corresponding to the action of $\mathbb{Q}(\widehat{C_r^*(\mathbb{Z} \times \mathbb{Z}_2)}, \clf)$, the fundamental unitary is 
         $$\begin{pmatrix}
a_1^{11} & a_2^{11}  & a_1^{12} & a_2^{12}\\
(a_2^{11})^\ast & (a_1^{11})^\ast & (a_2^{12})^\ast & (a_1^{12})^\ast\\
a_1^{21} & a_2^{21}  & a_1^{22} & a_2^{22}\\
(a_2^{21})^\ast & (a_1^{21})^\ast & (a_2^{22})^\ast & (a_1^{22})^\ast\\
\end{pmatrix}.$$

         We break the proof into a number of steps.  
         \paragraph{Step 1 :} 
         Using the condition  $\alpha(\gamma.\gamma^*)= \alpha(\gamma^\ast.\gamma)$, comparing the coefficients of $ \gamma^2, (\gamma^\ast)^2$ and $\gamma.\gamma^*$  on both sides we find that
         \begin{eqnarray}
         a_1^{11}(a_2^{11})^\ast &=& (a_2^{11})^\ast a_1^{11}, \label{eq 7.2} \\
          a_2^{11}(a_1^{11})^\ast &=& (a_1^{11})^\ast a_2^{11}, \label{eq 7.3} \\
          a_1^{11}(a_1^{11})^\ast + a_2^{11}(a_2^{11})^\ast &=& (a_2^{11})^\ast a_2^{11} + (a_1^{11})^\ast a_1^{11}. \label{eq 7.4} 
          \end{eqnarray}
            Applying the antipode on (\ref{eq 7.4}) we obtain
            \begin{equation}\label{eq 7.5}
            a_1^{11}(a_1^{11})^\ast + (a_2^{11})^\ast a_2^{11} = a_1^{11}(a_2^{11})^\ast + (a_1^{11})^\ast a_1^{11}.
            \end{equation}
            Thus one can conclude that both the elements $a_1^{11}$ and $a_2^{11}$ are normal using (\ref{eq 7.4}) and (\ref{eq 7.5}).  Hence the $C^*$-algebra $C^*\lbrace a_1^{11}, a_2^{11} \rbrace $ is commutative by (\ref{eq 7.2}), (\ref{eq 7.3}) and Proposition \ref{normal comp}.\\
             Applying the same argument replacing $\gamma$ by $\gamma^{\prime}$ one can deduce that $C^*\lbrace a_1^{22}, a_2^{22} \rbrace $  is commutative as well.\\
              Moreover, using $\alpha(\gamma.\gamma^*)= \alpha(\gamma^\ast.\gamma)
$, comparing the coefficients of $(\gamma^{\prime})^2,(\gamma^{\prime\ast})^2$ and $\gamma^{\prime}\gamma^{\prime\ast}$ on both sides, we have
\begin{eqnarray}
a_1^{12}(a_2^{12})^\ast &=& (a_2^{12})^\ast a_1^{12}, \label{eq 7.6}\\
 a_2^{12}(a_1^{12})^\ast &=& (a_1^{12})^\ast a_2^{12}, \label{eq 7.7}\\
 a_1^{12}(a_1^{12})^\ast + a_2^{12}(a_2^{12})^\ast &=& (a_2^{12})^\ast a_2^{12} + (a_1^{12})^\ast a_1^{12}. \label{eq 7.8}
 \end{eqnarray}
  Applying $\kappa$ on (\ref{eq 7.8}) we get
  \begin{equation} \label{eq 7.9}
 a_1^{21}(a_1^{21})^\ast + (a_2^{21})^\ast a_2^{21} = a_2^{21}(a_2^{21})^\ast + (a_1^{21})^\ast a_1^{21}.
 \end{equation}
 On the other hand, one find
 \begin{equation} \label{eq 7.10}
 a_1^{21}(a_1^{21})^\ast + a_2^{21}(a_2^{21})^\ast = (a_2^{21})^\ast a_2^{21} + (a_1^{21})^\ast a_1^{21},
 \end{equation}
 by comparing the coefficient of $\gamma.\gamma^*$ on both sides of $ \alpha(\gamma^{\prime})\alpha(\gamma^{\prime\ast})= \alpha(\gamma^{\prime\ast})\alpha(\gamma^{\prime})$. Now, by (\ref{eq 7.9}) and (\ref{eq 7.10}), both the elements $a_1^{21},a_2^{21}$ are normal. Using the antipode $a_1^{12},a_2^{12}$ are seen to 
be normal too. 
Thus, the $C^*$-algebra generated by $\lbrace a_1^{12},a_2^{12} \rbrace$ is commutative. Hence, $C^*\lbrace a_1^{21},a_2^{21} \rbrace$ is commutative as well, by applying the antipode.
       \paragraph{Step 2 :} 
        From the given fact $\alpha(\gamma.\gamma^{\prime})=0$, comparing the coefficients of $\gamma^2$ and $(\gamma^{\ast})^2$, we obtain
        \begin{equation}\label{eq 7.11}
         a_1^{11}a_1^{21}= a_2^{11}a_2^{21}=0.
         \end{equation}
         Applying the antipode on (\ref{eq 7.11}) and using the Proposition \ref{normal comp}  we have
         \begin{equation}
         a_1^{11}a_1^{12}= a_2^{11}a_2^{12}=0.
         \end{equation}
          Now, repeating the similar arguments using the relation $\alpha(\gamma.(\gamma^{\prime})^*)=0$ one can easily check that $a_1^{11}a_2^{12}= a_2^{11}a_1^{12}=0$. It can also be shown that $a_1^{22}a_1^{21}= a_1^{22}a_2^{21}=a_2^{22}a_1^{21}= a_2^{22}a_2^{21}=0$ by following the same line of arguments.
        
      \paragraph{Step 3 :}       
          Comparing the coefficient of $\gamma^2$ from the  relation $\alpha(\gamma.\gamma^\ast + \gamma^{\prime}.\gamma^{\prime\ast}) = 1 \ot 1$ we have
          \begin{equation}\label{eq 7.12}
          a_1^{11}(a_2^{11})^* + a_1^{21}(a_2^{21})^* = 0.
          \end{equation}
Multiplying by $a_2^{11}$ on the right side of the equation (\ref{eq 7.12}) one can get $a_1^{11}(a_2^{11})^*a_2^{11}=0$ as $(a_2^{21})^*a_2^{11}=0$. Hence $a_1^{11}(a_2^{11})^*=0$, which shows that  $a_1^{11}a_2^{11}=0$ by the Proposition \ref{normal comp} and also $a_1^{21}a_2^{21}=0$. Applying the same argument and comparing $(\gamma^\prime)^2$ from the relation  $\alpha(\gamma.\gamma^\ast + \gamma^{\prime}.\gamma^{\prime\ast}) = 1 \ot 1$ one can show that $a_1^{12}a_2^{12}= a_1^{22}a_2^{22}=0$. \\

Note that the underlying $C^*$-algebra of $\mathbb{Q}(\mathbb{Z})\star \mathbb{Q}(\mathbb{Z})$ is isomorphic to $(C(\mathbb{T}) \oplus C(\mathbb{T})) \star (C(\mathbb{T}) \oplus C(\mathbb{T}))$ which is the same as $C^*\lbrace \gamma_1,\gamma_2\rbrace \star C^*\lbrace \gamma_1^{\prime},\gamma_2^{\prime}\rbrace$, where $\lbrace \gamma_1,\gamma_2\rbrace$ and $\lbrace \gamma_1^{\prime},\gamma_2^{\prime}\rbrace$  satisfy the relations (\ref{eq 7.1}), (\ref{eq 71}) and (\ref{eq 72}). Finally, using the above steps we define the $C^*$-isomorphism from $(\mathbb{Q}(\mathbb{Z})\star \mathbb{Q}(\mathbb{Z})) \oplus (\mathbb{Q}(\mathbb{Z})\star \mathbb{Q}(\mathbb{Z}))$ to $\mathbb{Q}(\widehat{C_r^*(\mathbb{Z} \times \mathbb{Z}_2)}, \clf)$ by
      $$(\gamma_i,0)\mapsto a_i^{11}, (\gamma_i^{\prime},0)\mapsto a_i^{22}, (0,\gamma_i)\mapsto a_i^{21}, (0,\gamma_i^{\prime})\mapsto a_i^{12},$$
       for $i=1,2$. Indeed this is the doubling of $\mathbb{Q}(\mathbb{Z})\star \mathbb{Q}(\mathbb{Z})$ with respect to the order $2$ automorphism $\theta$ defined as in the statement of Theorem \ref{sec 7 th}. \ \qed\\  
         \brmrk
         As  $\mathbb{Q}(\mathbb{Z}) \star \mathbb{Q}(\mathbb{Z})$ is noncommutative, it is clear that the quantum isometry group of $C_r^*(\mathbb{Z} \times \mathbb{Z}_2)$ with the new generating set $\clf$ differs from the previous one, calculated in Theorem \ref{thm 3 for dir}  with group-like generating elements.
         Moreover, $\mathbb{Q}(\widehat{C_r^*(\mathbb{Z} \times \mathbb{Z}_2)}, \clf)$ can also be identified with the quantum group $K_2^+$ (for more details see Section $5$ of \cite{org filt}). We also give yet another description of $K_n^+$ in \cite{mandal 2}.  
          \ermrk

\end{document}